\documentclass{amsproc}[12pt]
\usepackage{amssymb}
\usepackage{amsfonts,amsmath,amssymb,amsthm}
\usepackage{latexsym,amscd}
\usepackage{amsbsy}
\usepackage{xy}
\xyoption{all}
\usepackage{bm}
\usepackage{graphicx,psfrag,subfigure}
\usepackage{amsfonts}

 \xyoption{all} \oddsidemargin0.6cm \evensidemargin0.6cm
\newcommand{\bcen}{\begin{center}}      \newcommand{\ecen}{\end{center}}

\def\Ext{\mbox{Ext}}
\def\rank{\mbox{rank}}
\def\dim{\mbox{dim}}
\def\HH{\mbox{HH}}
\def\Hom{\mbox{Hom}}
\def\Im{\mbox{Im}}

\def\Ker{\mbox{Ker}}
\def\wt{\widetilde}
\def\ch{\mbox{char}}
\def\mod{\mbox{mod}}

\def\bc{\begin{center}}
\def\ec{\end{center}}
\def\no{\noindent}
\def\hang{\hangindent\parindent}
\def\textindent#1{\indent\llap{\qquad #1\ \ \enspace}\ignorespaces}
\def\ref{\par\hang\textindent}

\begin{document}
\abovedisplayskip=6pt plus 1pt minus 1pt \belowdisplayskip=5pt plus
1pt minus 1pt
\thispagestyle{empty} \vspace*{-1.0truecm} \noindent
\vskip -8mm

\bc{\large\bf More counterexamples to Happel's question and
Snashall-Solberg's conjecture}
  \ec


\vskip 2mm \bc{\bf Xu Yunge$^{a,}$\footnote[1]{Corresponding
author:\\  e-mail: xuyunge@yahoo.com.cn (Xu Yunge),
zhangchao198658@yahoo.com.cn (Zhang Chao)} and  Zhang
Chao$^{a,b}$\\[5pt]
{\tiny  $^a$ School of Mathematics and Computer Science, Hubei
University, Wuhan 430062, PR China \\ $^b $ Academy of Mathematics
and Systems Science, Chinese Academy of Sciences, Beijing 100190, PR
China }}\ec

\vskip 1 mm

\noindent{\small {\small\bf Abstract} \ \ In this paper we provide
more counterexamples to Happel's question and Snashall-Solberg's
conjecture which generalize many counterexamples to these
conjectures studied in the literature. In particular, we show that a
family of $\mathbb{Z}_n\times \mathbb{Z}_n$-Galois covering algebras
of quantized exterior algebra $A_q$ in two variables answer
negatively to Happel's question, and meanwhile, the one-point
coextensions of $\mathbb{Z}_n$ and $\mathbb{Z}_n\times
\mathbb{Z}_m$-Galois covering algebras of $A_q$ negate the
Snashall-Solberg's conjecture.

\vspace{1mm}\baselineskip 12pt

\no{\small\bf Keywords} \ \ Hochschild cohomology ring; \ Happel's
question; \ Snashall-Solberg's conjecture;\ Koszul  dual; \ graded
center

\no{\small\bf MR(2000) Subject Classification}\ \ {\rm 16E40,
16G10}}

\section{ Introduction}

  Let
$\Lambda$ be a finite-dimensional $k$-algebra (associative with
identity) over a field $k$. Denote by $\Lambda^e$ the enveloping algebra of
$\Lambda$, i.e., the tensor product $\Lambda \otimes_k \Lambda^{op}$
of the algebra $\Lambda$ and its opposite $\Lambda^{op}$. Then by
Cartan-Eilenberg \cite{CE} the $i$-th Hochschild  homology and
cohomology groups of $\Lambda$ are identified with the $k$-spaces
 $$ \HH_i(\Lambda)=\mbox{Tor}_i^{\Lambda^e}(\Lambda, \Lambda),\quad
 \HH^i(\Lambda)=\Ext^i_{\Lambda^e}(\Lambda, \Lambda) $$ respectively.
 The Hochschild cohomology ring
 $\HH^{\ast}(\Lambda)=\bigoplus_{i=0}^{\infty}\HH^i(\Lambda)$
 has a so-called Gerstenhaber algebra structure under
 the cup product and the Gerstenhaber bracket \cite{Ger}. It is well known, as a
noncommutative analogy of differential forms and polyvector fields,
that Hochshild homology and cohomology of an associative
(noncommutative) algebra have been a starting point of
noncommutative geometry and play an important role due to the
classic  Hochschild-Kostant-Rosenberg theorem.

It is also well known that the homological properties of an algebra
are closely related to the properties of its Hochschild (co)homology
groups. For example, if a finite dimensional algebra over an
algebraically closed field has finite global dimension, then all its
higher Hochschild cohomology groups vanish.  The inverse is known as
Happel's question  and it has been shown that the conjecture does
not hold for the quantized exterior algebra $A_q=k\langle
x,y\rangle/(x^2, xy+qyx, y^2)$ (or more generally, for the quantized
complete intersection) when $q\in k^*=k\setminus\{0\}$ is not a root
of unity in \cite{Hap,BGMS,BE,Opp}. However, the homology version of
Happel's question comes to be known as ``Han's conjecture" and
remains still open \cite{Han}.

 Motivated by support variety for finitely
generated modules over group algebras defined by Carlson in
\cite{Car}, Snashall and Solberg developed support variety theory of
finitely generated modules over a finite-dimensional algebra in
\cite{SS}. They found that the finiteness condition of Hochschild
cohomology ring modulo nilpotence $\HH^{\ast}(\Lambda)/\mathcal{N}$
played an important role in support variety theory, where
$\mathcal{N}$ denotes the ideal of $\HH^{\ast}(\Lambda)$ generated
by all homogeneous nilpotent elements. Moreover, they also
conjectured that the Hochschild cohomology ring modulo nilpotence of
any finite-dimensional algebra is  a finitely generated  algebra,
and the conjecture was proved to be true for many classes of
algebras, such as finite-dimensional algebras of finite global
dimension\cite{Hap}, finite-dimensional monomial algebras\cite{GS,
GSS1}, finite-dimensional self-injective algebras of finite
representation type over an algebraically closed field\cite{GSS2},
any block of a group ring of a finite group\cite{Eve, Ven} and so
on. Until 2008, Xu F. provided the first counterexample to the
conjecture by studying the Hochschild cohomology ring modulo
nilpotence of a  seven-dimensional category algebra in the case of
$\ch k=2$ \cite{XuF}, which is isomorphic to a Koszul algebra
\cite{Sna}. Furthermore, it was proved that the Hochschild
cohomology ring modulo nilpotence of the above Koszul algebra as
well as its quantized algebra is not a finitely
 generated algebra irrespective of the characteristic of the base field $k$\cite{Sna,XZ}.

Let $\Lambda_q$ be the algebra introduced in the first paragraph of
the section 2, which arises from a formal deformation with
infinitesimal in $\HH^2(\Lambda)$ and occurs in the study of the
Drinfeld double of the generalized Taft algebras and of the
representation theory of $U_q(\textsf{sl}_2)$. Snashall and
Taillefer proved the Hochschild cohomology rings modulo nilpotence
of $\Lambda$ and $\Lambda_q$ are finitely generated commutative
algebras of Krull dimension two and hence Snashall-Solberg's
conjecture holds for this class of algebras\cite{ST0,ST}. However,
 Parker and Snashall showed in \cite{SP} that $\Lambda_q$ is an
infinite family of counterexamples to Happel's question  when
$\zeta=q_0q_1\cdots q_{m-1}$ is not a root of unity. Furthermore, we
prove that, for the algebra $\Gamma_q$ introduced in the first
paragraph of the section 2 which can be viewed as a one-point
coextension of $\Lambda_q$, $\HH^{\ast}(\Gamma_q)/\mathcal{N}$ is
not a finitely generated algebra if $\zeta=q_0q_1\cdots q_{m-1}$ is
a root of unity and thus provides an infinite family of
counterexamples to Snashall-Solberg's conjecture.

Note that, when $q_0= q_1= \cdots= q_{m-1}$, the algebra $\Lambda_q$
is just a $\mathbb{Z}_n$-Galois covering algebra of the quantized
exterior algebra $A_q$ \cite{HZ}, while the algebra $\Gamma_q$ can
be viewed as a one-point coextension of the $\mathbb{Z}_n$-covering
algebra $\Lambda_q$. So it seems that the following question arises
naturally: {\it if an algebra $A$ (for example, the quantized
exterior algebra $A_q$) answers negatively to Happel's question,
does it so for any finite Galois covering algebra $\wt{A}$ of $A$,
and meanwhile, will the one-point (co)extension of $\wt{A}$ provide
a family of counterexamples to the Snashall-Solberg's conjecture?}
Let $G$ be a finite group, $A$ a $G$-graded $k$-algebra, and
$\wt{A}$ the covering algebra with the Galois group $G$. It was
shown in \cite{MMM1,MMM2} that there is a ring monomorphism from
$\HH^{i}(\wt{A})$ to $\HH^{i}(A)$ for $i\geq  0$. As a consequence,
if $A$ is a counterexample to Happel's question, then so is
$\wt{A}$. Indeed, this is the case for the $\mathbb{Z}_2$-graded
quantized exterior algebra $A_q$ and its Galois covering algebra
with Galois group $\mathbb{Z}_2$ (even more generally,
$\mathbb{Z}_n$) \cite{XT,SP}. However, if $A$ is only a $k$-algebra
(unnecessarily $G$-graded), then there is only a monomorphism from
$\HH^{i}(\wt{A})^G$ to $\HH^{i}(A)$ for $i\geq  0$, and the explicit
descriptions of these maps for $i=0, 1$ are provided in \cite{GHS}.

In this paper, we first employ Snashall and Taillefer's strategy in
\cite{ST} to consider the structure of the Hochschild cohomology
rings modulo nilpotence $\HH^{\ast}(\Lambda_q)/\mathcal{N}$ of the
algebras $\Lambda_q$  by computing the graded center of its Koszul
dual $E(\Lambda_q)$ in the section 2. As a consequence, we show that
they are not finitely generated as algebras when $\zeta$ is a root
of unity, and thus provide more counterexamples to
Snashall-Solberg's conjecture, which include and generalize all the
counterexamples studied in \cite{XuF,Sna,XZ}. Next,  we consider a
family of algebras $\Lambda_q^{m,n}$ as well as their one-point
coextensions $\Gamma_q^{m,n}$, where
$q=(q_{00},q_{01},\cdots,q_{n-1, m-1})\in (k^*)^{mn}$. In the case
that $q_{ij}=q_{00}$ for all $i,j$, $\Lambda_q^{m,n}$ is just a
covering algebra
 of the quantized exterior algebra $A_q$ with the Galois group
 $\mathbb{Z}_n\times \mathbb{Z}_m$.
We determine the structure of Hochschild cohomology ring of
$\Lambda_q^{n,n}$ and show that $\Lambda_q^{n,n}$ answers negatively
to Happel's question when $\xi=\prod_{i,j=0}^{n-1}q_{ij}\in k^*$ is
not a root of unity in the section 3, and meanwhile,
$\Gamma_q^{m,n}$ provides an infinity family of counterexamples to
Snashall-Solberg's conjecture when
$\eta=\prod_{i=0}^{n-1}\prod_{j=0}^{m-1}q_{ij}\in k^*$ is a root of
unit in the section 4 as expected.

\medskip
\section{ Graded center of $E(\Gamma_q)$}

Throughout this section we always assume that  $\Lambda=kQ/I$ is a
class of selfinjective Koszul algebras, where the quiver $Q$ is of
the form in the left hand side of Figure 1 below,  and the ideal $I$
is generated by the set $R=\{a_{i}a_{i+1},{b_{i-1}}{b_{i-2}},
a_{i}b_{i}-b_{i-1}a_{i-1}\mid 0 \leq i\leq m-1\}$. Let $\Lambda_q$,
$q=(q_0, q_1, \cdots, q_{m-1})\in (k^*)^{m}$, be their socle
deformations (i.e. $\Lambda_q$ are selfinjective with
$\Lambda_q/\mbox{soc}(\Lambda_q)\cong \Lambda/\mbox{soc}(\Lambda)$),
see also \cite{ST0,ST}. Throughout we always assume that all the
subscripts of arrows are identified with their residues modulo $m$.

Let $\Gamma_q=kQ'/I'_q$, where $Q'$ is the finite quiver with $m+1$
vertices $\{0,1,\ldots,m-1\}\cup\{-1\}$ and $3m$ arrows pictured in
Figure 1 as well, and $I'_q$ is the ideal generated by
$R'=\{a_{i}a_{i+1},{b_{i+1}}{b_{i}}, q_{i}a_{i}b_{i}-b_{i-1}a_{i-1},
a_ic_{i+1}\mid 0 \leq i\leq m-1,a_m=a_0\}$ and $q=(q_0, q_1, \cdots,
q_{m-1})\in (k^*)^m$. In the case of $m=1$,  $\Gamma_q$ is
isomorphic to the quantized Koszul algebra studied in \cite{XZ} (in
which the ``commutative" relation is $ab+qba$) and used to provide a
family of counterexamples to Snashall-Solberg's conjecture.
Throughout the section we assume $m\geq 2$.

\begin{figure}
\includegraphics[height=150pt,width=180pt]{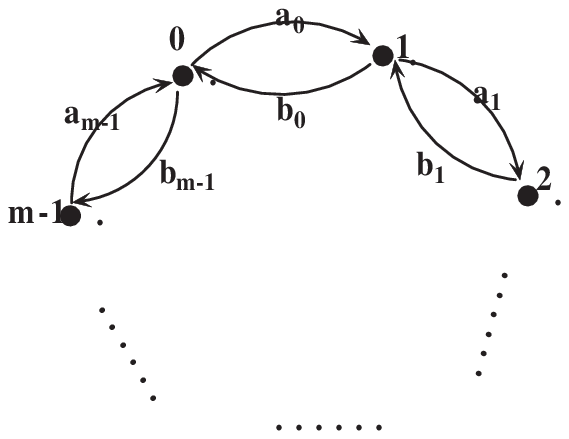} \qquad
\includegraphics[height=150pt,width=180pt]{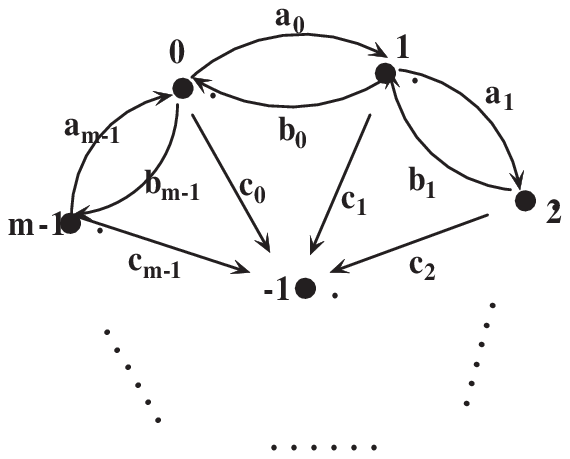}
\caption{The quivers $Q$ and
  $Q'$}
\label{fig:eg-1}
\end{figure}

 Denote by $e_i$ the trivial
path in $kQ'$ and write the composition of arrows from left to
right. Note that the left length lexicographic order by choosing
$e_0<\cdots<
e_{m-1}<e_{-1}<a_0<\cdots<a_{m-1}<b_0<\cdots<b_{m-1}<c_0<\cdots<c_{m-1}$
provides an admissible order on $kQ'$ and thus the set $R'$ is just
a (noncommutative) quadratic reduced Gr\"obner basis of
$I'_q$\cite{DMR}. So $\Gamma_q$ is Koszul by \cite{GH}.

Recall that the Ext-algebra $E(\Gamma_q)$ is just the Koszul dual of
$\Gamma_q$. Thus $E(\Gamma_q)=kQ^{op}/{I'}_q^{\perp}$, where
${I'}_q^{\perp}=\langle q_i^{-1}(a_ib_i)^o+(b_{i-1}a_{i-1})^o,\,
(b_ic_i)^o\rangle$ and $x^o\in Q^o$ denotes the opposite arrow of
the arrow $x$ in $Q$. Moreover, any left $kQ^{op}$-module can be
viewed as a right $kQ$-module, so we may consider $E(\Gamma_q)$ as
the quotient of $kQ$  modulo the ideal generated by
$q_i^{-1}a_ib_i+b_{i-1}a_{i-1}, b_ic_i$ for $i=0, 1, \cdots, m-1$.

In a similar way to \cite{ST}, we denote by $\gamma^s_i$ and
$\delta_i^t$ the paths $a_ia_{i+1}\cdots a_{i+s-1}$ and
$b_{i+t-1}\cdots b_{i+1}b_i$ respectively. Unless otherwise
specified, we do not distinguish a path with its image in
$E(\Gamma_q)$. Thus any typical monomial in $E(\Gamma_q)$ has the
form $\gamma^s_i\delta^t_j$ for some integers $s, t$ and $0\leq i,j
\leq m-1$. The algebra $E(\Gamma_q)$ is a bigrading algebra graded
with  the lengths of paths and with the degree induced  by choosing
the degree of $e_i$, $a_j$, $c_j$ and $b_j$ to be $0,1,1,-1$
respectively. Thus any monomial element $\gamma^s_i\delta^t_j$ has
the length $s+t$ and degree $s-t$. we denote by $|z|$ the length of
a length-homogeneous element $z$ in $E(\Gamma_q)$. Denote by
$Z_{gr}(E(\Gamma_q))$  the graded center of $E(\Gamma_q)$.

It is easy to see that $z\in Z_{gr}(E(\Gamma_q))$ if and only if $z$
satisfies the following conditions:

(1) $e_jz=ze_j$, for $j=-1,0, 1, \cdots, m-1$;

(2) $a_jz=(-1)^{|z|}za_j$, for $0\leq j\leq m-1$;

 (3) $b_jz=(-1)^{|z|}zb_j$, for any $0\leq j\leq m-1$;

(4) $c_jz=(-1)^{|z|}zc_j$, for any $0\leq j\leq m-1$.

{\bf Lemma 2.1.} If a homogeneous element $z\in
Z_{gr}(E(\Gamma_q))$, then $z\in k$, or $z$ has the form
 $$z=\sum_{i=0}^{m-1}u_i\gamma^{s_0}_i\delta^{t_0}_i, \quad u_i\in k$$ with $s_0\equiv t_0(\mod ~
 m)$, $t_0\geq 1$, and for $0\leq j\leq m-1$,
$$u_{j+1}=(-1)^{s_0}(q_{j+1}\cdots
 q_{j+t_0})^{-1}u_j=(-1)^{t_0}(q_{j+1}\cdots q_{j+s_0})^{-1}u_j.\eqno{(2\mbox{-}1)}$$

 {\bf Proof.}
By (1), we can write $z=\sum_{i=-1}^{m-1}e_ize_i$. Note that for any
$0\leq i\leq
 m-1$, a typical monomial in $e_iE(\Gamma_q)e_i$ has the form $\gamma^s_i\delta^t_i$, where
 $s, t\geq 0$, and $s\equiv t(\mod \, m)$. In particular,
 $e_{-1}E(\Gamma_q)e_{-1}=e_{-1}$. Moreover, $Z_{gr}(E(\Gamma_q))$ is generated by the elements
 which are both length homogeneous and degree homogeneous.  Therefore, if the length of
 $z$ is $0$, then $z=\sum_{i=-1}^{m-1}d_ie_i$, where $d_i\in k$; otherwise,
 $z$ has the form $\sum_{i=0}^{m-1}u_i\gamma^{s_i}_i\delta^{t_i}_i,$
 where $u_i\in k$, $s_i, t_i\geq 0$, $s_i\equiv t_i(\mod \, m)$.
 The degree homogeneity implies that
 $s_i-t_i=s_0-t_0$ and the length homogeneity implies that $s_i+t_i=s_0+t_0>0$,
 and hence we have $s_i=s_0$ and $t_i=t_0$ for $i=0, 1, \cdots,
 m-1$. So
 $$z=\sum_{i=0}^{m-1}u_i\gamma^{s_0}_i\delta^{t_0}_i.$$ Here we also take the subscripts modulo
 $m$ (in particular, $u_0=u_m$).

 We next consider the condition (2).  If the length of $z=\sum_{i=-1}^{m-1}d_ie_i$ is zero, then
 $d_{j+1}a_j=a_j\sum_{i=-1}^{m-1}d_ie_i=(\sum_{i=-1}^{m-1}d_ie_i)a_j=d_ja_j$,
 and we have $z=d_0\sum_{i=0}^{m-1}e_i+d_{-1}e_{-1}$. On the other
 hand, if the length of $z$ is not zero, we have
 $$a_jz=u_{j+1}a_j\gamma^{s_0}_{j+1}\delta^{t_0}_{j+1}=
 u_{j+1}\gamma^{s_0+1}_j\delta^{t_0}_{j+1}$$ and
 $$za_j=u_j\gamma^{s_0}_j\delta^{t_0}_ja_j=
 (-1)^{t_0}u_j(q_{j+1}\cdots
 q_{j+t_0})^{-1}\gamma^{s_0+1}_j\delta^{t_0}_{j+1}.$$ The condition
 (2) implies that  $u_{j+1}=(-1)^{s_0}(q_{j+1}\cdots
 q_{j+t_0})^{-1}u_j\neq 0$ and similarly, the condition (3) implies that  $u_{j+1}=(-1)^{t_0}(q_{j+1}\cdots
 q_{j+s_0})^{-1}u_j\neq 0$ for $0\leq j\leq m-1$.

 By the condition (4), we know that if the length of $z$ is
 zero, then $d_{-1}c_j=c_j(d_0\sum_{i=0}^{m-1}e_i+d_{-1}e_{-1})
 =(d_0\sum_{i=0}^{m-1}e_i+d_{-1}e_{-1})c_j=d_0c_j$, so $d_{-1}=d_0$, and thus
 $z=d_0(\sum_{i=-1}^{m-1}e_i)=d_0$. Otherwise,
 $z=\sum_{i=0}^{m-1}u_i\gamma^{s_0}_i\delta^{t_0}_i$ satisfies
 $0=c_jz=(-1)^{s_0+t_0}zc_j=(-1)^{s_0+t_0}\sum_{i=0}^{m-1}u_i\gamma^{s_0}_i\delta^{t_0}_ic_j$
for all $0\leq j\leq m-1$ in $E(\Gamma_q)$,  which forces that
$t_0\geq 1$ by the definition of ${I'}_q^{\perp}$. We complete the
proof of the lemma. \hfill$\square$

\medskip
{\bf Remark.} Comparing with the result in \cite{ST}, we have
$Z_{gr}(E(\Gamma_q))\setminus
k=\{z=\sum_{i=0}^{m-1}u_i\gamma^{s_0}_i\delta^{t_0}_i\in
Z_{gr}(E(\Lambda_q))\mid t_0\geq 1\}$. Using the formula (2-1)
recursively, one can obtain that
$$u_i=(-1)^{is_0}\prod_{k=1}^{i}(q_k\cdots
q_{k+t_0-1})^{-1}u_0=(-1)^{it_0}\prod_{k=1}^{i}(q_k\cdots
q_{k+s_0-1})^{-1}u_0,$$ for $i=1,2,\dots, m-1$. In particular,
  $$u_0=u_m=(-1)^{ms_0}(q_0\cdots q_{t_0-1})^{-1}(q_1\cdots
 q_{t_0})^{-1}\cdots (q_{m-1}\cdots q_{m-2+t_0})^{-1}u_0.$$ Since $u_0\neq
 0$, we have $(q_0\cdots q_{t_0-1})^{-1}(q_1\cdots
 q_{t_0})^{-1}\cdots (q_{m-1}\cdots q_{m-2+t_0})^{-1}=(-1)^{ms_0}$.
 Let $\zeta=q_0q_1\cdots q_{m-1}$, then
 $\zeta^{t_0}=(-1)^{ms_0}$. Similarly,  $\zeta^{s_0}=(-1)^{mt_0}$.

 \medskip

 {\bf Proposition 2.2.}~ If $\zeta$ is not a root of unity, then
 $Z_{gr}(E(\Gamma_q))=k$.

 {\bf Proof.}~For any element $z\in Z_{gr}(E(\Gamma_q))$, if the length of $z$ is not zero,
 then $\zeta^{t_0}=(-1)^{ms_0}$ and
 $\zeta^{s_0}=(-1)^{mt_0}$. Since $\zeta$ is not a root of unity,
 then $s_0=t_0=0$, this yields a contradiction. Thus the length of $z$
 is zero. By Lemma 2.1 we have $z\in k$. On the other hand,
 it is evident that $k\subseteq Z_{gr}(E(\Gamma_q))$, therefore,
 $Z_{gr}(E(\Gamma_q))=k$ as desired. \hfill$\square$

 \medskip

Now we assume that $\zeta$ is a primitive $d$-th root of unity, that
is, $d\geq 1$ is the minimal integer such that $\zeta^d=1$.  The
proof of the following  proposition is similar to that of
\cite[Prop.2.4, 2.5]{ST} and hence we omit all the details and leave
only the sketch of the proof.

 \medskip

{\bf Proposition 2.3.}~Suppose that $\zeta$ is a primitive $d$-th
root of unity. Then
$$Z_{gr}(E(\Gamma_q))\cong\left\{\begin{array}{ll}
 (k[x, y, w]/\langle w^{2m}-\epsilon_d xy\rangle)_{x^*},
 &\mbox{if }m \mbox{ is odd, char}k\ne 2, \mbox{ and }d\equiv 2(\mod~4);\\
 (k[x, y, w]/\langle w^{m}-\epsilon_d xy\rangle)_{x^*},
 &\mbox{otherwise.}\end{array}\right. $$ where $(k[x, y, w]/\langle
w^p-\epsilon_d xy\rangle)_{x^*}$ denotes the subalgebra of $k[x, y,
w]/\langle w^p-\epsilon_d xy\rangle$ that does not contain the
subspace spanned by the $x^i$ for $i=1,2,\cdots$, and
$$\epsilon_d=\left\{\begin{array}{lll}
\prod_{l=1}^{m-1}\prod_{k=1}^{ld}(-1)^{md/2}(q_k\cdots
q_{k+d-1})^{-1},& \mbox{if } m\mbox{ is even or char}k=2;\\
\prod_{l=1}^{m-1}\prod_{k=1}^{ld}(q_k\cdots
q_{k+d-1})^{-1},& \mbox{if }m\mbox{ is odd, char}k\ne 2,\mbox{ and } ~d\equiv 0(\mod ~4);\\
\prod_{l=1}^{2m-1}\prod_{k=1}^{ld/2}(q_k\cdots
q_{k+d/2-1})^{-1},& \mbox{if }m\mbox{ is odd, char}k\ne 2,\mbox{ and } d\equiv 2(\mod ~4);\\
\prod_{l=1}^{m-1}\prod_{k=1}^{2ld}(q_k\cdots q_{k+2d-1})^{-1},&
\mbox{if }m\mbox{ is odd, char}k\ne 2,\mbox{ and } d \mbox{ is odd}.
\end{array}\right.$$

{\bf Proof.}\footnote{For the referee's convenience, we leave the
complete proof of the proposition in the appendix.} Clearly,
$k\subseteq Z_{gr}(E(\Gamma_q))$. Note that
$Z_{gr}(E(\Gamma_q))\setminus
k=\{z=\sum_{i=0}^{m-1}u_i\gamma^{s_0}_i\delta^{t_0}_i\in
Z_{gr}(E(\Lambda_q))\mid u_i\in k, t_0\geq 1\}$ by the remark above.

~{\bf Case 1.} $m$ is even or char$k=2$.  With a similar but lengthy
analysis as in \cite[Prop.2.4]{ST}, any homogeneous element $z\in
Z_{gr}(E(\Gamma_q))\setminus k$ can be written as
$$z=u_0\Big(\sum_{i=0}^{m-1}(-1)^{is}\prod_{k=1}^{i}(q_k\cdots
q_{k+s-1})^{-1}\gamma^{s}_i
\delta^{s}_i\Big)\Big(\sum_{i=0}^{m-1}\gamma^{dm}_i\Big)^{\alpha
}\Big(\sum_{i=0}^{m-1}\delta^{dm}_i\Big)^{\alpha+h},\eqno{(2\mbox{-}2)}$$
where $s=ld$ for $0\leq l\leq m-1$. Set
$w=\sum_{i=0}^{m-1}(-1)^{id}\prod_{k=1}^{i}(q_k\cdots
q_{k+d-1})^{-1}\gamma^{d}_i \delta^{d}_i$,
$x=\sum_{i=0}^{m-1}\gamma^{md}_i$, $y=\sum_{i=0}^{m-1}\delta^{md}_i$
and
$\epsilon_d=(-1)^{md/2}\prod_{l=1}^{m-1}\prod_{k=1}^{ld}(q_k\cdots
q_{k+d-1})^{-1}$. Then $w^m=\epsilon_d xy$. Moreover, by the formula
(2-2),  any homogeneous element
$z=\sum_{i=0}^{m-1}u_i\gamma^{s_0}_i\delta^{t_0}_i\in
Z_{gr}(E(\Gamma_q))\setminus k$ can be written as a scalar multiple
of $x^iy^jw^l$ with $j+l>0$ since $t_0\geq 1$. The condition $j+l>0$
implies that $x^i, i\geq 1$, does not belong to
$Z_{gr}(E(\Gamma_q))$. With the same argument as in \cite[Lemma
2.3]{ST},  the elements $x, y, w$ don't have additional relation
except $w^m=\epsilon_d xy$ in $Z_{gr}(E(\Gamma_q))$. So
$Z_{gr}(E(\Gamma_q))\cong  (k[x, y, w]/\langle w^{m}-\epsilon_d
xy\rangle)_{x^*}$.

{\bf Case 2.} $m$ is odd and char$k\ne 2$.

(i)~If $d$ is odd, writing $t_0=\alpha dm+t, s_0=\beta dm+s$, then
any homogeneous element $z\in Z_{gr}(E(\Gamma_q))\setminus k$ can be
written as
\begin{eqnarray*}z&=&
u_0\sum_{i=0}^{m-1}\prod_{k=1}^{i}(q_k\cdots
q_{k+dl/2-1})^{-1}\gamma^{dl/2}_i\delta^{dl/2}_i(\sum_{i=0}^{m-1}\gamma^{2dm}_i)^{\alpha/2}
(\sum_{i=0}^{m-1}\delta^{2dm}_i)^{\beta/2}, \quad\mbox{or}\\
z&=&u_0\sum_{i=0}^{m-1}\prod_{k=1}^{i}(q_k\cdots
q_{k+dl/2-1})^{-1}\gamma^{d(l/2+m)}_i\delta^{d(l/2+m)}_i(\sum_{i=0}^{m-1}\gamma^{2dm}_i)^{(\alpha-1)/2}
(\sum_{i=0}^{m-1}\delta^{2dm}_i)^{(\beta-1)/2}
\end{eqnarray*}
depending on whether $\alpha$ is even or odd.  Set
$w=\sum_{i=0}^{m-1}\prod_{k=1}^{i}(-1)^{i}(q_k\cdots
q_{k+2d-1})^{-1}\gamma^{2d}_i\delta^{2d}_i$,
$x=\sum_{i=0}^{m-1}\gamma^{2md}_i$,
$y=\sum_{i=0}^{m-1}\delta^{2md}_i$ and
$\epsilon_d=\prod_{l=1}^{m-1}\prod_{k=1}^{2dl}(q_k\cdots
q_{k+2d-1})^{-1}$. Then $w^m=\epsilon_d xy$. Moreover, any
homogeneous element $z\in Z_{gr}(E(\Gamma_q))\setminus k$  can be
written as a scalar multiple of $x^{i}y^{j}w^l$ with $j+l>0$. And
there is no additional relation in $Z_{gr}(E(\Gamma_q))$ except
$w^{m}=\epsilon_d xy$.

 (ii)~If $d$ is even and  $d\equiv 0(\mod ~4)$,
 then  any homogeneous element $z\in Z_{gr}(E(\Gamma_q))\setminus k$
 can be written as
$$z=u_0\sum_{i=0}^{m-1}\prod_{k=1}^{i}(q_k\cdots
q_{k+dl/2-1})^{-1}\gamma^{dl/2}_i\delta^{dl/2}_i(\sum_{i=0}^{m-1}\gamma^{dm}_i)^{\alpha}
(\sum_{i=0}^{m-1}\delta^{dm}_i)^{\beta}.
$$
Set $w=\sum_{i=0}^{m-1}\prod_{k=1}^{i}(q_k\cdots
q_{k+d-1})^{-1}\gamma^{d}_i\delta^{d}_i$,
$x=\sum_{i=0}^{m-1}\gamma^{md}_i$, $y=\sum_{i=0}^{m-1}\delta^{md}_i$
and $\epsilon_d=\prod_{l=1}^{m-1}\prod_{k=1}^{dl}$ $(q_k\cdots
q_{k+d-1})^{-1}$. Then $w^m=\epsilon_d xy$. And we can write
homogeneous element $z\in Z_{gr}(E(\Gamma_q))$ which is not in $k$
as a scalar multiple of $x^{\alpha}y^{\beta}w^{l/2}$ with $\beta
+l/2>0$. In particular, any scalar multiple of $x^i$ does not lie in
$Z_{gr}(E(\Gamma_q))$, for $i=1, 2, \cdots$. Also, there is no
additional relation in $Z_{gr}(E(\Gamma_q))$ except
$w^{m}=\epsilon_d xy$.

(iii)~ If $d$ is even and $d\equiv 2(\mod ~ 4)$, then  we can write
any homogeneous element $z\in Z_{gr}(E(\Gamma_q))$ that is not in
$k$ as
$$z =u_0\sum_{i=0}^{m-1}\prod_{k=1}^{i}(-1)^{il}(q_k\cdots
q_{k+dl/2-1})^{-1}\gamma^{dl/2}_i\delta^{dl/2}_i(\sum_{i=0}^{m-1}\gamma^{dm}_i)^{\alpha}
(\sum_{i=0}^{m-1}\delta^{dm}_i)^{\beta}.
$$
Set $w=\sum_{i=0}^{m-1}\prod_{k=1}^{i}(-1)^{i}(q_k\cdots
q_{k+d/2-1})^{-1}\gamma^{d/2}_i\delta^{d/2}_i$,
$x=\sum_{i=0}^{m-1}\gamma^{md}_i$, $y=\sum_{i=0}^{m-1}\delta^{md}_i$
and $\epsilon_d=\prod_{l=1}^{2m-1}\prod_{k=1}^{dl/2}(q_k\cdots
q_{k+d/2-1})^{-1}$, then $w^{2m}=\epsilon_d xy$. And we can write
homogeneous element $z\in Z_{gr}(E(\Gamma_q))$ which is not in $k$
as a scalar multiple of $x^{\alpha}y^{\beta}w^{l}$ with $\beta+l>0$,
which implies that any scalar multiple of $x^i$ does not lie in
$Z_{gr}(E(\Gamma_q))$, for $i=1, 2, \cdots$. Again, there is no
additional relation in $Z_{gr}(E(\Gamma_q))$ except
$w^{2m}=\epsilon_d xy$. \hfill$\square$

\medskip

By \cite{SS, BGSS}, we know that
$\HH^{\ast}(\Gamma_q)/\mathcal{N}\cong
Z_{gr}(E(\Gamma_q))/\mathcal{N}_Z$, where $\mathcal{N}_Z$ denotes
the ideal of $Z_{gr}(E(\Gamma_q))$ generated by all nilpotent
elements. It follows directly from the above two propositions that
$\mathcal{N}_Z=0$. As a result, we have, in fact, characterized the
structure of $\HH^{\ast}(\Gamma_q)/\mathcal{N}$ and provided more
counterexamples to Snashall-Solberg's conjecture by the following
theorem.

{\bf Theorem 2.4.}~ Let $q=(q_0, q_1, \cdots, q_{m-1})\in
(k^{\ast})^{m}$, and $\zeta=q_0q_1\cdots q_{m-1}$. If $\zeta$ is a
not root of unity, then $\HH^{\ast}(\Gamma_q)/\mathcal{N}\cong k$;
If $\zeta$ is a root of unity, then
$\HH^{\ast}(\Gamma_q)/\mathcal{N}\cong Z_{gr}(E(\Gamma_q))$ is not
finitely generated as algebra.

{\bf Proof.} From the proposition 2.3, we know that if $\zeta$ is a
root of unity, then $\HH^{\ast}(\Gamma_q)/\mathcal{N}\cong (k[x, y,
w]/\langle w^p-\epsilon xy\rangle)_{x^*}$. Note that $x^iy$ lies in
$(k[x, y, w]/\langle w^p-\epsilon xy\rangle)_{x^*}$  but $x^i$ does
not, for $i=1, 2, \cdots$, then  $x^iy$ can not be generated by the
elements of lower degree in $(k[x, y, w]/\langle w^p-\epsilon
xy\rangle)_{x^*}$, and thus $\HH^{\ast}(\Gamma_q)/\mathcal{N}$ is
not finitely generated algebra when $\zeta$ is a root of unity.
\hfill$\square$

\section{The Hochschild cohomology ring of $\Lambda^{m,n}_q$}

Throughout this section, we assume $\Lambda^{m,n}_q=k\bar{Q}/I_q$,
where $\bar{Q}$ is a torus-like finite quiver which has $mn$
vertices $\{(i,j)\mid i\in \mathbb{Z}_n, j\in \mathbb{Z}_m \}$, and
$2mn$ arrows: $\{a_{ij}: (i,j)\rightarrow (i, j+1)\}\cup \{b_{ij}:
(i,j)\rightarrow (i+1, j)\}$ pictured as in Figure 2, and
$I_q=\langle a_{ij}a_{i,j+1},\; b_{ij}b_{i+1, j},\;
a_{ij}b_{i,j+1}+q_{ij}b_{ij}a_{i+1,j}\mid i\in \mathbb{Z}_n, j\in
\mathbb{Z}_m\rangle$, $q_{ij}\in k^{*}$. Denote by $e_{ij}$ the
idempotent of $\Lambda^{m,n}_q$ at the vertex $(i,j)$. Note that
$\Lambda^{m,n}_q$ is the $\mathbb{Z}_n\times \mathbb{Z}_m$-Galois
covering algebra  of the quantized exterior algebra $A_q$ if
$q_{ij}=q_{00}$ for $i\in \mathbb{Z}_n, j\in \mathbb{Z}_m$.

For $x\in \{e, a, b\}$, define $x_{ij}<x_{pl}$ if and only if $i<p$
or $i=p$ but $j<l$; and set $e_{i_1j_1}<a_{i_2j_2}<b_{i_3j_3}$.
 Then the
length-left-lexicographic order provides an admissible order for
$k\bar{Q}$, and $R=\{a_{ij}a_{i,j+1}, b_{ij}b_{i+1, j},
a_{ij}b_{i,j+1}+q_{ij}b_{ij}a_{i+1,j}\}$ forms a noncommutative
quadratic reduced Gr{\"o}ber basis of the ideal $I_q=\langle
a_{ij}a_{i,j+1}, b_{ij}b_{i+1, j},
a_{ij}b_{i,j+1}+q_{ij}b_{ij}a_{i+1,j}\rangle$, thus
$\Lambda^{m,n}_q$ is a Koszul algebra\cite{DMR, GH}.

In this section, we first construct a minimal projective bimodule
resolution of $\Lambda^{m,n}_q$, and then determine the structure of
Hochschild cohomology ring of $\Lambda^{m,n}_q$ when $m=n$ and
$\xi=\prod_{i,j=0}^{n-1}q_{ij}$ is not a root of unity, and thus
provide another family of counterexamples to Happel's conjecture.
For the convenience of notations, we denote by $\Lambda^n_q$ the
algebra $\Lambda^{n, n}_q$ unless otherwise specified in this
section. If $n=1$, then $\Lambda_q^1$ is just the quantized exterior
algebra $A_q$; if $n=2$, the Hochschild homology and cohomology of
$\Lambda_q^2$ have been considered in \cite{HX} and the
$k$-dimension of $\HH^*(\Lambda_q^2)$ is $4$ in the case that $q$ is
not a root of unity. From now on we assume $n\geq  3$ in this
section.

Let \begin{eqnarray*}&&g^0=\{g^0_{0,i,j}=e_{ij}\};\\
&&g^1=\{g^1_{0,i,j}=a_{ij}, \; g^1_{1,i,j}=b_{ij}\};\\
&&g^2=\{g^2_{0,i,j}=a_{ij}a_{i,j+1}, \;
g^2_{1,i,j}=a_{ij}b_{i,j+1}+q_{ij}b_{ij}a_{i+1,j},\;
g^2_{2,i,j}=b_{ij}b_{i+1,j}\}.\end{eqnarray*} Moreover, we set
$g^l_{-1,i,j}=0=g^l_{l+1,i,j}$, and when $l\geq 3$,
$g^l=\{g^l_{pij}\mid 0\leq p\leq l\}$, where
$$g^l_{pij}=a_{ij}g^{l-1}_{p,i,j+1}+q_{ij}q_{i,j+1} \cdots
q_{i, j+l-p-1}b_{ij}g^{l-1}_{p-1, i+1, j}.~~\eqno{(3\mbox{-}1)}$$

Define $P_l=\Lambda_q^n\otimes_{E} kg^l \otimes_{E} \Lambda_q^n$,
where $\Lambda_q^n=E\oplus \textrm{r}$ and $E\cong
\Lambda_q^n/\textrm{r}\cong k\times k\times \cdots \times k$. We
denote $\otimes_E$ by $\otimes$ for legibility. Set
$\tilde{g}^l_{pij}=1\otimes g^l_{pij}\otimes 1$ for $0\leq p\leq l$,
$l=0, 1, 2, \cdots$ and define $d_l: \, P_l\rightarrow P_{l-1}$ for
$l\geq 1$ as follows
\begin{eqnarray*}d_l(\tilde{g}^l_{pij})
&=&a_{ij}\tilde{g}^{l-1}_{p,i,j+1}+q_{ij}q_{i,j+1} \cdots q_{i,
j+l-p-1}b_{ij}\tilde{g}^{l-1}_{p-1, i+1,
j}+(-1)^l\tilde{g}^{l-1}_{p-1,i,j}b_{i+p-1, j+l-p}\\
&&+(-1)^lq_{i, j+l-p-1} \cdots q_{i+p-1, j+l-p-1}\tilde{g}^{l-1}_{p,
i, j}a_{i+p, j+l-p-1}.
\end{eqnarray*}

\medskip

\begin{figure}
\includegraphics[height=200pt,width=300pt]{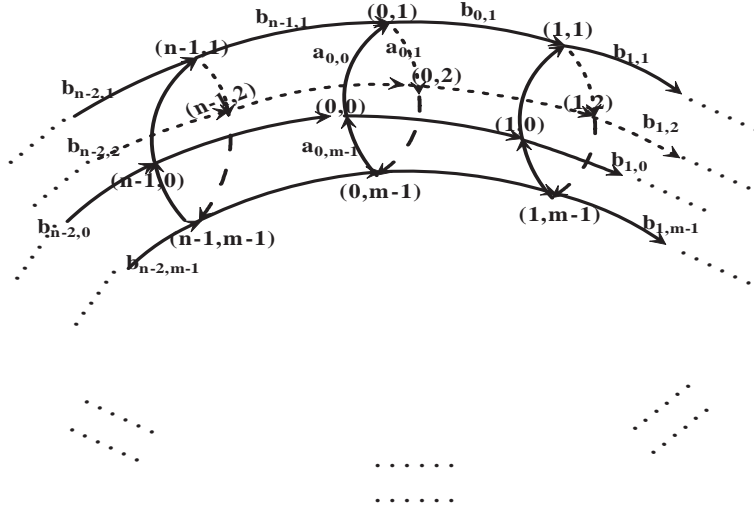}
\caption{The quiver $\bar{Q}$} \label{fig:eg-1}
\end{figure}

{\bf Lemma 3.1.} The complex $(\mathbb{P},d)$
$$
 \cdots\rightarrow P_{l+1}
\stackrel{d_{l+1}}{\longrightarrow} P_l
\stackrel{d_l}{\longrightarrow} \cdots
\stackrel{d_3}{\longrightarrow} P_2\stackrel{d_2}{\longrightarrow}
P_1 \stackrel{d_1}{\longrightarrow} P_0 \rightarrow 0
$$ is a minimal projective bimodule resolution of $\Lambda_q^n$.

{\bf Proof.} Let $X=\mbox{span}\{ a_{ij}, b_{ij}\mid i, j\in
\mathbb{Z}_n\}$. Since $\Lambda_q^n$ is  Koszul,  it suffices to
prove that the set $g^l$ forms a $k$-basis of
$K_l:=\cap_{s+t=l-2}X^sRX^t$ by \cite[Sec.9]{BK}, where, by abuse of
notation, $R$ stands for the $k$-space spanned by the set
$\{a_{ij}a_{i,j+1}, b_{ij}b_{i+1, j},
a_{ij}b_{i,j+1}+q_{ij}b_{ij}a_{i+1,j}\mid i, j\in \mathbb{Z}_n\}$.

We will first show that $g^l\subseteq K_l$ by induction on $l$. It
is clear when $l=2$. Assume that it holds for $l-1$. It is not
difficult but lengthy to verify that
$$g^l_{pij}=g^{l-1}_{p-1,i,j}b_{i+p-1, j+l-p}+q_{i, j+l-p-1} \cdots q_{i+p-1,
j+l-p-1}g^{l-1}_{p, i, j}a_{i+p, j+l-p-1}, \; 0\leq p\leq l.
$$
Then by the formula (3-1) and the above formula, we have
$g^l\subseteq XK_{l-1} \cap K_{l-1}X\subseteq K_{l}$.

On the other hand,  each element $g^l_{pij}$ in $g^l$ contains
exactly $p$ many $b$-class arrows, and hence the elements in $g^l$
are linearly independent. Moreover, the Koszul dual of $\Lambda_q^n$
is just the quadratic dual $k\bar{Q}^{op}/I_q^{\bot}$, where
$I_q^{\bot}=\langle
(a_{ij}b_{i,j+1})^{o}-q^{-1}_{ij}(b_{ij}a_{i+1,j})^{o}\rangle$, so
the Betti numbers of a minimal projective bimodule resolution of
$\Lambda_q^n$ are $\{b_l=(l+1)n^2\}$, and thus $\dim K_l=(l+1)n^2$,
which coincides with the number of elements in $g^l$.

The differential $d$ is obtained from \cite{GHMS} directly. The
proof is completed. \hfill$\square$

\medskip

In order to compute the Hochschild cohomology  of $\Lambda_q^n$ when
$\xi$ is not a root of unity, we first recall some notations from
\cite{Cib90}. We say a path $\alpha$ is {\it uniform} if there exist
$(i, j), (i', j')\in \bar{Q}_0$, such that $\alpha=e_{ij}\alpha
e_{i'j'}$. Two paths $\alpha$ and $\beta$ are said to be {\it
parallel}, and denoted by $\alpha//\beta$, provided that they have
the same source and target. If $X$ and $Y$ are  sets of some uniform
paths in $\bar{Q}$, then $X// Y:=\{(\alpha, \beta)\in X\times Y\mid
\alpha//\beta \}$ and we denote by $k(X // Y)$ the $k$-vector space
with the set $X
// Y$ as basis.

Applying the functor $\Hom_{(\Lambda_q^n)^e}(-,
 \Lambda_q^n)$ to the minimal projective bimodule resolution
 $(\mathbb{P},d)$ of $\Lambda_q^n$, we obtain the Hochschild
 cochain complex  $C^{\ast}(\mathbb{P})$:
$$\quad 0\rightarrow
\Hom_{(\Lambda_q^n)^e}(P_0, \Lambda_q^n)
\stackrel{d^1}{\longrightarrow}
\cdots \stackrel{d^n}{\longrightarrow}\Hom_{(\Lambda_q^n)^e}(P_n, \Lambda_q^n)\\
\stackrel{d^{n+1}}{\longrightarrow} \Hom_{(\Lambda_q^n)^e}(P_{n+1} ,
\Lambda_q^n) \stackrel{d^{n+2}}{\longrightarrow} \cdots,
$$ where $d^i:=\Hom_{(\Lambda_q^n)^e}(d_i, \Lambda_q^n), i=0,1,2,\cdots$.

Let $\mathcal{B}=\{e_{ij}, a_{ij}, b_{ij}, a_{ij}b_{i, j+1}\mid i,j\in \mathbb{Z}_n\}$ be a
$k$-basis of $\Lambda_q^n$. Thanks to the isomorphism in
\cite{Cib90}, that is, $k(g^l// \mathcal {B})\stackrel{\phi}{\cong}
\Hom_{(\Lambda_q^n)^e}(P_l, \Lambda_q^n)$, where $\phi: (g^l_{pij},
x)\mapsto f_{(g^l_{pij}, x)}$, $x\in \mathcal {B}$, and
$f_{(g^l_{pij}, x)}(1\otimes g^l_{p^{'}i^{'}j^{'}}\otimes
1)=\delta_{pij,p^{'}i^{'}j^{'}}x. $ Here $\delta_{pij,p^{'}i^{'}j^{'}}$ denotes
the Kronecker sign, that is,  $\delta_{pij,p^{'}i^{'}j^{'}}=1$
if $(p,i,j)=(p',i',j')$ (i.e. $p=p^{'}, i=i^{'}, j=j^{'}$) and 0 otherwise. Under the
isomorphism  the complex $(C^{\ast}(\mathbb{P}), d^{\ast})$ changes
into
$$(M^{\bullet}, \delta^{\bullet})= 0\rightarrow
k(g^0//\mathcal {B}) \stackrel{\delta^1}{\longrightarrow} \cdots
\stackrel{\delta^l}{\longrightarrow}k(g^l//\mathcal {B})
\stackrel{\delta^{l+1}}{\longrightarrow} k(g^{l+1}//\mathcal {B})
\stackrel{\delta^{l+2}}{\longrightarrow} \cdots,
$$  where
$$\begin{array}{ll} & \delta^l(g^{l-1}_{pij},
x)=\phi^{-1}d^l\phi(g^{l-1}_{pij}, x)\\ = & (g^l_{p,i,j-1},
a_{i,j-1}x)+q_{i-1,j}\cdots q_{i-1, j+l-p-2}(g^l_{p+1,i-1,j},
b_{i-1,j}x)
 +(-1)^l(g^{l}_{p+1, i, j}, xb_{i+p, j+l-p-1})\\
 & +(-1)^lq_{i, j+l-p-1}
\cdots q_{i+p-1, j+l-p-1}(g^{l}_{p, i, j}, xa_{i+p, j+l-p-1}).
\end{array}$$

By definition, we know that
$\HH^l(\Lambda_q^n)=\Ker\delta^{l+1}/\Im\delta^l$, thus
\begin{eqnarray*}\dim_k\HH^l(\Lambda_q^n)&=&\dim_k\Ker\delta^{l+1}-\dim_k\Im\delta^l\\
&=&\dim_kM^l-\dim_k\Im\delta^{l+1}-\dim_k\Im\delta^l.
\end{eqnarray*}
Since the set $\mathcal{B}=\{e_{ij}, a_{ij}, b_{ij}, a_{ij}b_{i,
j+1}\}$ is a $k$-basis of $\Lambda_q^n$, the elements in
$(g^l//\mathcal{B})$ has the form of $(g^l_{pij}, x)$ with $x\in
\mathcal{B}$. Note that $l$ stands for the length of $g^l_{pij}$ and
$p$ describes the number of $b$-class arrows appearing in each
monomial of $g^l_{pij}$, and $g^l_{pij}$ is uniform with the source
$(i, j)$ and the target $(i+p, j+l-p)$. Thus $(g^l_{pi'j'}//e_{ij})$
implies $(i', j')=(i,j)$ and $l=l_0n$, $p=un$ for some integers
$u,l_0$ with $0\leq u\leq l_0$. Similarly, the elements  in
$g^{\bullet}$ parallel to $a_{ij}, b_{ij}, a_{ij}b_{i, j+1}$ have
the form of $g^{l_0n+1}_{un,i,j}, g^{l_0n+1}_{un+1,i,j},
g^{l_0n+2}_{un+1,i,j}$ respectively, where $0\leq u\leq l_0$.
Therefore,
$$\dim_k M^l= \dim_k (\Gamma_l//\mathcal{B})=\left\{\begin{array}{llll}
(l_0+1)n^2, \;&\mbox{if} \; l=l_0n;\\
2(l_0+1)n^2, \;&\mbox{if} \; l=l_0n+1;\\
 (l_0+1)n^2, \;&\mbox{if} \;l=l_0n+2;\\
0, \;&\mbox{otherwise.}\end{array}\right. \eqno{(3\mbox{-}2)}$$
Since
$M^{l_0n+2}=k(g^{l_0n+2}//\mathcal{B})=k\{(g^{l_0n+2}_{un+1,i,j},
a_{ij}b_{i, j+1})\mid 0\leq u\leq l_0, i,j\in \mathbb{Z}_n\}$, we
have $\delta^{l_0n+3}=0$  by the definition of $\delta^{\bullet}$.
Also, for $3<i\leq n$, $\delta^{l_0n+i}=0$ since $M^{l_0n+i-1}=0$.
Thus the complex $(M^{\bullet}, \delta^{\bullet})$ has the forms of
$$0\rightarrow
M^0 \stackrel{\delta^1}{\longrightarrow} M^1
\stackrel{\delta^2}{\longrightarrow} M^2
\stackrel{0}{\longrightarrow}\cdots
\stackrel{0}{\longrightarrow}M^{l_0n}
\stackrel{\delta^{l_0n+1}}{\longrightarrow} M^{l_0n+1}
\stackrel{\delta^{l_0n+2}}{\longrightarrow} M^{l_0n+2}
\stackrel{0}{\longrightarrow}\cdots,
$$ where $M^l=k(g^l//\mathcal{B})$.
So it suffices to consider $\dim_k\Im\delta^{l_0n+1}$ and
$\dim_k\Im\delta^{l_0n+2}$.

The order $<$ on $\mathcal{B}$ induces an order on
$(g^l//\mathcal{B})$ as follows: $(g^l_{pij}, x)\prec(g^l_{p'i'j'},
x')$ if and only if $p<p'$, or $p=p'$ but $x<x'$. By abuse of
notation, we denote  still by $\delta^{l}$ the matrix of the
differential $\delta^{l}$ under the ordered bases above. Then by the
description of $\delta^{l}$, $\delta^{l_0n+1}$ and $\delta^{l_0n+2}$
have the following form respectively:
$$
\delta^{l_0n+1}=\left(
  \begin{array}{cccccccc}
    A_0 &  & &    \\
      B_0 &  &  &    \\
     & A_1 &  &   \\
     & B_1 &  &    \\
    &  & \ddots&    \\
     &  & \ddots  &    \\
     &  & & A_{l_0}  \\
     &  &  &  B_{l_0}
  \end{array}
\right); \qquad  \delta^{l_0n+2}=\left(
  \begin{array}{cccccccc}
    C_0 & D_0 & &   &  &  & &  \\
     & &  C_1 &  D_1   &  &&  &   \\
     &  & &  & \ddots & \ddots  & &  \\
     &  &  &   & &  & C_{l_0} &  D_{l_0}
  \end{array}
\right),
$$ where $A_i=\mbox{diag}\{A_{i0},A_{i1}, \cdots A_{i, n-1}\}$,
and if we set $r_i=\prod_{j=0}^{n-1}q_{ij}$,
$c_j=\prod_{i=0}^{n-1}q_{ij}$, then
$$A_{ij}=\left(
  \begin{array}{cccccccc}
    (-1)^lc_0^i &1  & &    \\
     & (-1)^lc_1^i & \ddots &    \\
     &  & \ddots & 1  \\
     1&  &  &  (-1)^lc_{n-1}^i
  \end{array}
\right)_{n\times n};$$
$$B_i=\left(
  \begin{array}{cccc}
    (-1)^lI_n &r_0^{l_0-i}I_n  & &     \\
     & (-1)^lI_n &\ddots  &     \\
     &  & \ddots & r_{n-2}^{l_0-i}I_n  \\
    r_{n-1}^{l_0-i}I_n &  &  &  (-1)^lI_n
  \end{array}
\right)_{n^2\times n^2};
$$
$$C_i=\left(
  \begin{array}{cccc}
    (-1)^lI_n &-r_0^{l_0-i}I_n  & &     \\
     & (-1)^lI_n &\ddots  &     \\
     &  & \ddots & -r_{n-2}^{l_0-i}I_n  \\
    -r_{n-1}^{l_0-i}I_n &  &  & (-1)^lI_n
  \end{array}
\right)_{n^2\times n^2},
$$ where $I_n$ denotes the identity matrix of size $n\times n$; $D_i=\mbox{diag}\{D_{i0}, D_{i1}, \cdots,
D_{i,n-1}\}$, and $$D_{ij}=\left(
  \begin{array}{cccccccc}
    (-1)^{l+1}c_0^i &1  & &    \\
     & (-1)^{l+1}c_1^i & \ddots &    \\
     &  & \ddots & 1 \\
     1&  &  &  (-1)^{l+1}c_{n-1}^i
  \end{array}
\right)_{n\times n}.$$

{\bf Lemma 3.2.}~ Suppose that $n\geq  3$ and $\xi$ is not a root of
unity. Then  $$\dim_k\Im\delta^{l_0n+1}=\dim_k\Im\delta^{l_0n+2}=\left\{\begin{array}{ll}n^2-1, &\mbox{if  } l_0=0;\\
(l_0+1)n^2, &\mbox{otherwise. }\end{array}\right.$$

{\bf Proof.}~We first consider the matrix $\delta^{l_0n+1}$. Since,
for $0<i\leq l_0$,
$\det(A_i)=(\det(A_{i0}))^n=((-1)^{nl}\xi^i+(-1)^{n+1})^n\neq 0$ by
the condition that $\xi$ is not a root of unity, the last $l_0n^2$
columns of $\delta^{l_0n+1}$  are linearly independent. We assert
that $$\rank \left(\begin{array}{cc}
  B_0 \\ \hline    A_0 \end{array}
\right)=\left\{\begin{array}{ll} n^2,& \mbox{if }l_0>0;\\ n^2-1, &
\mbox{if }l_0=0. \end{array}\right.$$ Indeed,  by adding
$(-1)^{l+1}r_i^{l_0}$-multiple of the $(i+1)$-th block-column of
$B_0$ to the $(i+2)$-th block-column of $B_0$ in turn for $i=0, 1,
\ldots, n-2$, we obtain
$$B_0\longrightarrow \left(\begin{array}{cccc}
    (-1)^lI_n & & &    \\
      &   \ddots &&  \\
     &  &(-1)^l I_n &     \\
 r_{n-1}^{l_0}I_n & \cdots & (-1)^{(l+1)(n-2)}(r_{n-1}r_{1}\cdots r_{n-3})^{l_0}I_n &  (-1)^{(l+1)(n-1)}\xi^{l_0}I_n+(-1)^lI_n
  \end{array}
\right).$$ Thus
$\det(B_0)=\big((-1)^{(l+1)(n-1)}\xi^{l_0}+(-1)^{l}\big)^n$, which
is nonzero if $l_0> 0$. Thus $\rank \left(\begin{array}{cc}
  B_0 \\ \hline    A_0 \end{array}
\right)= n^2$ in the case when  $l_0>0$. If $l_0=0$, by adding the
 $(i+1)$-th block-column  to the
$(i+2)$-th block-column  in turn for $i=0, 1, \ldots, n-2$, we
obtain
$$\left(\begin{array}{cc}
   B_0 \\ \hline
    A_0
  \end{array}
\right) \longrightarrow\left(\begin{array}{cccccccc}
    -I_n & & &    \\
      &   \ddots &&  \\
     &  &-I_n &     \\
   I_n & \cdots & I_n &  0 \\ \hline
   A_{00}& & & A_{00}   \\
      &   \ddots && \vdots  \\
     &  &A_{0, n-2} &  A_{0, n-2}   \\
  &  &  &  A_{0, n-1}
  \end{array}
\right).$$ Since $A_{00}=A_{01}=\cdots=A_{0,n-1}$ and $\rank\,
A_{00}=n-1$, we have $\rank \delta^{l_0n+1}=\rank\,\delta^{1}=n^2-1$
in the case when $l_0=0$ as desired. Therefore,
$$\dim_k\Im\delta^{l_0n+1}=\rank \,
\delta^{l_0n+1}=l_0n^2+\rank \left(\begin{array}{cc}
  B_0 \\ \hline    A_0 \end{array}
\right)=\left\{\begin{array}{ll}n^2-1, &\mbox{if  } l_0=0;\\
(l_0+1)n^2, &\mbox{otherwise. }\end{array}\right.$$ We complete the
proof of the first part of this lemma.

Next, we consider the rank of $\delta^{l_0n+2}$. With a similar
argument as for $\delta^{l_0n+1}$,
$\det(D_i)=(\det(D_{i0}))^n=((-1)^{n(l+1)}\xi^i+(-1)^{n+1})^n\neq 0$
for $0<i\leq l_0$ since $\xi$ is not a root of unity. Therefore, the
last $l_0n^2$ rows of $\delta^{l_0n+2}$ are linearly independent and
it suffices to consider the rank of $(D_0|C_0)$. We claim that
$$\rank \left(\begin{array}{c|c}
  D_0 &    C_0 \end{array}
\right)=\left\{\begin{array}{ll} n^2,& \mbox{if }l_0>0;\\ n^2-1, &
\mbox{if }l_0=0. \end{array}\right.$$ In fact,  by adding
$(-1)^{l}r_i^{l_0}$-multiple of the $(i+2)$-th block-row of $C_0$ to
the $(i+1)$-th block-row of $C_0$ in turn for $i=n-2, n-1, \ldots,
0$, we obtain
$$C_0\longrightarrow \left(\begin{array}{cccc}
     -(-1)^{l(n-1)}\xi^{l_0}I_n+(-1)^lI_n & & &    \\
     -(-1)^{l(n-2)}(r_{n-1}\cdots r_{1})^{l_0}I_n  &   (-1)^l I_n &&  \\
    \vdots &  & \ddots&     \\
 -r_{n-1}^{l_0}I_n & &  &(-1)^lI_n
  \end{array}
\right).$$ So
$\det(C_0)=\big(-(-1)^{l(n-1)}\xi^{l_0}+(-1)^{l}\big)^n\neq 0$,  if
$l_0> 0$. Thus $\rank (D_0|C_0)= n^2$ in the case when  $l_0>0$. If
$l_0=0$, then, by adding the
 $(i+1)$-th block-row  to the
$(i+2)$-th block-row  in turn for $i=0, 1, \ldots, n-2$, we obtain
$$(D_0|C_0) \longrightarrow\left(\begin{array}{cccc|cccc}
    D_{00} & & & &I_n&&&   \\
      &   \ddots && &\vdots & \ddots &&  \\
   &  &D_{0,n-2}& &I_n& &I_n&     \\
    D_{00} & \cdots & D_{0,n-2} &  D_{0,n-1} &0&\cdots &0& 0
  \end{array}
\right).$$ Since $D_{00}=D_{01}=\cdots=D_{0,n-1}$ and $\rank\,
D_{00}=n-1$, we have $\rank \delta^{l_0n+2}=\rank\,\delta^{2}=n^2-1$
in the case when $l_0=0$, which proves our claim. Therefore,
$$\dim_k\Im\delta^{l_0n+2}=\rank \,
\delta^{l_0n+2}=l_0n^2+\rank (D_0|C_0)=\left\{\begin{array}{ll}n^2-1, &\mbox{if  } l_0=0;\\
(l_0+1)n^2, &\mbox{otherwise. }\end{array}\right.$$\hfill$\square$

\medskip

With the help of Lemma 3.2, we immediately have the following
theorem.

{\bf Theorem 3.3.}~ If $n\geq  3$ and $\xi$ is not a root of
unity, then we have $$\dim_k\HH^l(\Lambda_q^n)=\left\{\begin{array}{llll}1, &\mbox{if  } l=0 \mbox{ or  } 2;\\
2, &\mbox{if  } l=1;\\
0, &\mbox{otherwise. }\end{array}\right.$$ Thus $\HH^*(\Lambda_q^n)$
is a finite dimensional algebra of dimension $4$.

{\bf Proof.}~ It follows directly from Lemma 3.2 and the formula
$$\dim_k\HH^l(\Lambda_q^n)=\dim_kM^l-\dim_k\Im\delta^{l+1}-\dim_k\Im\delta^l.$$ \hfill$\square$

\medskip
{\bf Remark.} Note that our result still holds true for $n=2$ (cf.
\cite{HX}). Moreover, it also shows that, when $\xi$ is not a root
of unity, $\Lambda_q^n$ provides a family of counterexamples to
Happel's question as expected.

{\bf Corollary 3.4.}~ If $\xi$ is not a root of unity, then
$\HH^*(\Lambda_q^n)\cong \wedge(u,v)$, the exterior algebra.

{\bf Proof.}~ For legibility, we do not distinguish the parallel
path in $M^l$ with its image in $\HH^l(\Lambda^n_q)$. Moreover, it
is straightforward to calculate that
$\HH^0(\Lambda_q^n)=\mbox{span}\{\sum_{i,j}(g^0_{0ij},
e_{ij})\}\cong k$,
$\HH^1(\Lambda_q^n)=\mbox{span}\{\sum_{i,j}(g^1_{0ij}, a_{ij}),
\sum_{i,j}(g^1_{1ij}, b_{ij})\}$,
$\HH^2(\Lambda_q^n)=\mbox{span}\{\sum_{i,j}(g^2_{1ij},
a_{ij}b_{i,j+1})\}$. Under the isomorphism $\phi: k(g^l// \mathcal
{B})\rightarrow \Hom_{(\Lambda_q^n)^e}(P_l, \Lambda_q^n)$, we have
$f^1_a=\sum_{i,j}f_{(g^1_{0ij}, a_{ij})}$ and
$f^1_b=\sum_{i,j}f_{(g^1_{1ij}, b_{ij})}$ also form a $k$-basis of
$\HH^1(\Lambda_q^n)$, and $f^2_{ab}=\sum_{i,j}f_{(g^2_{1ij},
a_{ij}b_{i,j+1})}$ a $k$-basis of $\HH^2(\Lambda_q^n)$. We define
bimodule maps
$$\psi_0: P_1\rightarrow P_0, \quad \left\{\begin{array}{ll} \tilde{g}^1_{0,i,j}\mapsto & 0,\\
\tilde{g}^1_{1,i,j}\mapsto &
b_{ij}\tilde{g}^0_{0,i+1,j};\end{array}\right.$$
$$\psi_1: P_2\rightarrow P_1, \quad \left\{\begin{array}{ll} \tilde{g}^2_{0,i,j}\mapsto & 0,\\
\tilde{g}^2_{1,i,j}\mapsto & -q_{ij}
b_{ij}\tilde{g}^1_{0,i+1,j};\\
\tilde{g}^2_{2,i,j}\mapsto & -b_{ij}
\tilde{g}^1_{1,i+1,j}\end{array}\right.$$ Now it is easy to check
that the following diagram is commutative:
$$\xymatrix{
P_2\ar[r]^{d_2}\ar[d]_{\psi_1}& P_1\ar[d]_{\psi_0}\ar[rd]^{f^1_b}
&\\ P_1\ar[r]^{d_1}\ar[rd]_{f^1_a} & P_0\ar[r]^{\mu}& \Lambda^n_q\\
&\Lambda^n_q&}$$ where $\mu$ is the multiplication. Thus the
composition $f^1_a\psi_1: P_2\rightarrow \Lambda_q^n$ is just the
Yoneda product $f_a^1* f_b^1$ in $\HH^2(\Lambda_q^n)$, which is
$f^2_{ab}$ and thus is nonzero in $\HH^2(\Lambda_q^n)$. By the
graded commutativity of $\HH^*(\Lambda_q^n)$, we have $f_a^1*
f_b^1=- f_b^1* f_a^1$, and $f_a^1* f_a^1=0= f_b^1* f_b^1$ when
char$k\ne 2$. These still hold by a direct calculation when
char$k=2$. Denote $u=f^1_a$, $v=f^1_b$ for simplicity. So
$\HH^*(\Lambda_q^n)\cong \wedge(u,v)$.

\bigskip

\section{The graded center of $E(\Gamma^{m,n}_q)$}

Let $\Gamma^{m,n}_q=k\tilde{Q}/\tilde{I}_q$, where $\tilde{Q}$ is a
wheel-like finite quiver with $mn+1$ vertices:  $\{(i,j)\mid i\in
\mathbb{Z}_n, j\in \mathbb{Z}_m \}\cup \{-1\}$ , and $3mn$ arrows:
$\{a_{ij}: (i,j)\rightarrow (i, j+1)\}\cup \{b_{ij}:
(i,j)\rightarrow (i+1, j)\} \cup \{c_{ij}: (i,j)\rightarrow
-1\}$(see Figure 3), and $\tilde{I}_q=\langle a_{ij}a_{i,j+1},\;
b_{ij}b_{i+1, j}, \;a_{ij}c_{i,j+1},\;
a_{ij}b_{i,j+1}+q_{ij}b_{ij}a_{i+1,j}\rangle$, $q_{ij}\in k^{*}$. In
fact, the algebra $\Gamma_q^{m,n}$ can be regarded as a one-point
coextension of the algebra $\Lambda_q^{m,n}$ defined in the previous
section. Throughout this section, we assume that
$\eta=\prod_{i=0}^{n-1}\prod_{j=0}^{m-1}q_{ij}$, and denote by
$e_{ij}$ the idempotent of $\Lambda^{m,n}_q$ at the vertex $(i,j)$
and by $e_{-1}$ the idempotent at the vertex $-1$. In this section,
we will describe the graded center of $E(\Gamma^{m,n}_q)$ by
applying Snashall and Taillefer's method in \cite{ST} to the algebra
$\Gamma_q^{m,n}$.

In a similar way to the previous sections, we can show that
$\Gamma_q^{m,n}$ is a Koszul algebra. Moreover, its Koszul dual
$E(\Gamma_q^{m,n})=k\tilde{Q}^{op}/\tilde{I}_q^{\perp}$, where
$\tilde{I}_q^{\perp}=\langle (b_{ij}c_{i+1,j})^{o},\;
(a_{ij}b_{i,j+1})^{o}-q_{ij}^{-1}(b_{ij}a_{i+1,j})^{o}\rangle$ and
$x^{o}$ denotes the arrow in $\tilde{Q}^{op}$ corresponding to $x$
in $\tilde{Q}$. Moreover, $E(\Gamma_q^{m,n})$ can be viewed as a
quotient of $k\tilde{Q}$ modulo the ideal generated by
$b_{ij}c_{i+1,j}, a_{ij}b_{i,j+1}-q_{ij}^{-1}b_{ij}a_{i+1,j}$ for
$i\in \mathbb{Z}_n, j\in \mathbb{Z}_m$. Denote still by
$Z_{gr}(E(\Gamma_q^{m,n}))$  the graded center of
$E(\Gamma_q^{m,n})$.

In this section, we do not differentiate the path in $k\tilde{Q}$
and its image in $E(\Gamma_q^{m,n})$. Since $e_{ij}z=ze_{ij}$ and
$e_{-1}z=ze_{-1}$ for any $z\in Z_{gr}(E(\Gamma_q^{m,n}))$, we can
write
$z=\sum\limits_{i=0}^{n-1}\sum\limits_{j=0}^{m-1}e_{ij}ze_{ij}+e_{-1}ze_{-1}$.
Let $\alpha_{ij}$ and $\beta_{ij}$ denote the path
$a_{ij}a_{i,j+1}\cdots a_{i,j+m-1}$ and $b_{ij}b_{i+1,j}\cdots
b_{i+n-1,j}$ respectively.  Using the relation
$a_{ij}b_{i,j+1}=q_{ij}^{-1}b_{ij}a_{i+1,j}$ for $i\in \mathbb{Z}_n,
j\in \mathbb{Z}_m$ repeatedly, an element $z$ satisfying
$z=e_{ij}ze_{ij}$ can be written as the form
$z=u_{ij}\alpha_{ij}^{s_{ij}}\beta_{ij}^{t_{ij}}$ for some
$u_{ij}\in k$. Moreover, $e_{-1}ze_{-1}=e_{-1}$.

\begin{figure}
\includegraphics[height=200pt,width=300pt]{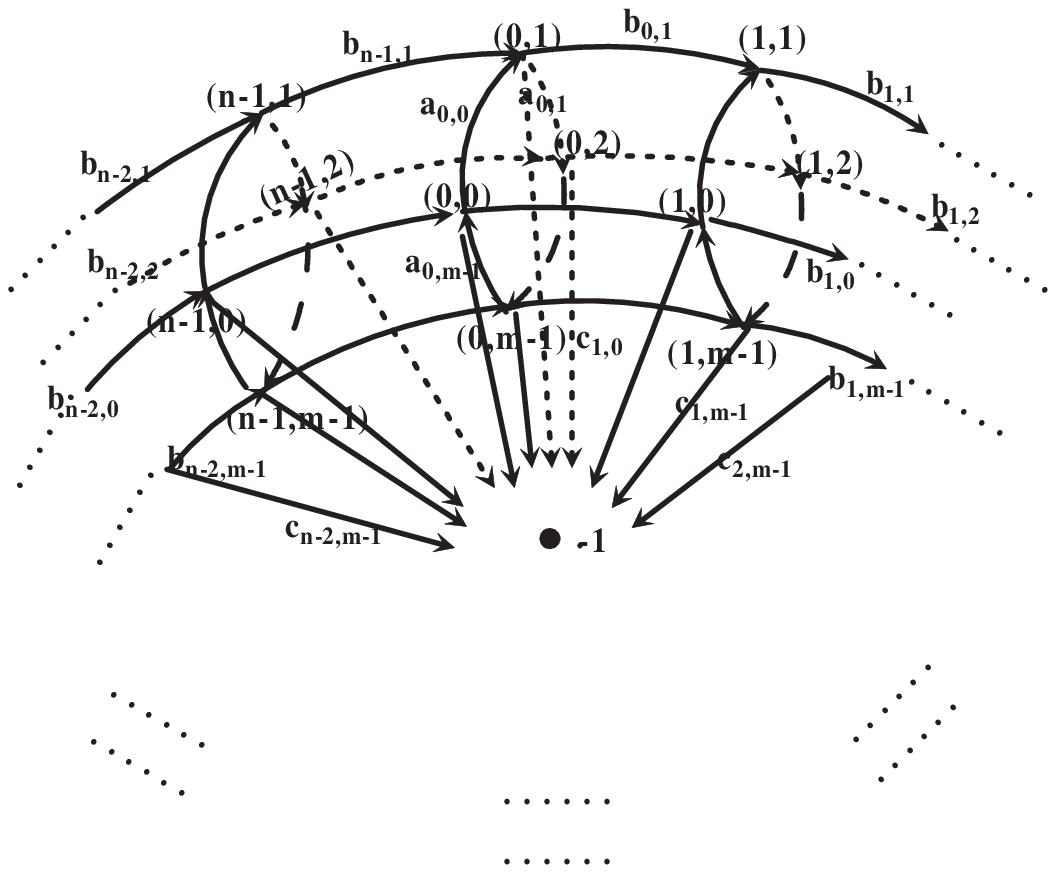}
\caption{The quiver $\tilde{Q}$} \label{fig:eg-1}
\end{figure}

Noting that $Z_{gr}(E(\Gamma_q^{m,n}))$ can be generated by some
elements which are length homogeneous, we denote by $|z|$ the length
of such an element $z$ and $z$ must satisfy the following additional
conditions:

(1) $a_{ij}z=(-1)^{|z|}za_{ij}$, for $i\in \mathbb{Z}_n, j\in
\mathbb{Z}_m$;

 (2) $b_{ij}z=(-1)^{|z|}zb_{ij}$, for any $i\in \mathbb{Z}_n,
j\in \mathbb{Z}_m$;

(3) $c_{ij}z=(-1)^{|z|}zc_{ij}$, for any $i\in \mathbb{Z}_n, j\in
\mathbb{Z}_m$.

\medskip

{\bf Lemma 4.1.}~  For any homogeneous element $z\in
Z_{gr}(E(\Gamma_q^{m,n}))$, we have $z\in k$ or $z$ can be written
as
$$z=\sum\limits_{i=0}^{n-1}\sum\limits_{j=0}^{m-1}u_{ij}\alpha_{ij}^{s_{0}}\beta_{ij}^{t_{0}}$$
with
$u_{ij}=(-1)^{(i+j)(ms_0+nt_0)}(\prod_{l=0}^{j-1}\prod_{p=0}^{n-1}q_{pl}^{t_0})
(\prod_{p=0}^{i-1}\prod_{l=0}^{m-1}q_{pl}^{-s_0})u_{00}\in k^*$ and
$t_0\geq 1$. Moreover,
$$\eta^{t_0}=(-1)^{m(ms_0+nt_0)},
\eta^{s_0}=(-1)^{n(ms_0+nt_0)}.$$

{\bf Proof.}~ We consider the condition (1). If $|z|=0$, then
$z=\sum_{i=0}^{n-1}\sum_{j=0}^{m-1}u_{ij}e_{ij}+u_{-1}e_{-1}$ with
$u_{ij}, u_{-1}\in k$, and $a_{ij}z=(-1)^{|z|}za_{ij}$ implies that
$u_{ij}=u_{i,j+1}$ for $i\in \mathbb{Z}_n, j\in \mathbb{Z}_m$. If
$|z|\ne 0$, then $z$ has the form
$z=\sum_{i=0}^{n-1}\sum_{j=0}^{m-1}u_{ij}\alpha_{ij}^{s_{ij}}\beta_{ij}^{t_{ij}}$
with $ms_{ij}+nt_{ij}=ms_{00}+nt_{00}$. Moreover,
$$a_{ij}z=a_{ij}\cdot
u_{i,j+1}\alpha_{i,j+1}^{s_{i,j+1}}\beta_{i,j+1}^{t_{i,j+1}}
=(q_{ij}\cdots
q_{i+n-1,j})^{-t_{i,j+1}}u_{i,j+1}\alpha_{ij}^{s_{i,j+1}}\beta_{ij}^{t_{i,j+1}}a_{ij},$$
and $za_{ij}=u_{ij}\alpha_{ij}^{s_{ij}}\beta_{ij}^{t_{ij}}a_{ij}$.
Thus the equality $a_{ij}z=(-1)^{|z|}za_{ij}$ implies that
$s_{ij}=s_{i,j+1}$, $t_{ij}=t_{i,j+1}$ and
$u_{i,j+1}=(-1)^{ms_{00}+nt_{00}}(q_{ij}\cdots
q_{i+n-1,j})^{t_{i,j+1}}u_{ij}$. Recursively, we have
$u_{0,0}=u_{0,m}=(-1)^{m(ms_{00}+nt_{00})}\eta^{t_{00}}u_{00}$.

Similarly, the condition (2) implies that if $|z|=0$, then
$u_{i+1,j}=u_{ij}$. Thus we have
$z=u_{00}\sum_{i=0}^{n-1}\sum_{j=0}^{m-1}e_{ij}+u_{-1}e_{-1}$ with
$u_{00}, u_{-1}\in k$. Moreover, if $|z|\ne 0$, then
$s_{i+1,j}=s_{ij}$, $t_{i+1,j}=t_{ij}$ and
$(-1)^{ms_{00}+nt_{00}}(q_{ij}\cdots
q_{i,j+m-1})^{-s_{ij}}u_{ij}=u_{i+1,j}$. Moreover, we  have that
$u_{0,0}=u_{n,0}=(-1)^{n(ms_{00}+nt_{00})}\eta^{-s_{00}}u_{00}$
recursively. For legibility of notations, we denote $s_{00}$ and
$t_{00}$ by $s_0$ and $t_0$ respectively. So, taking the condition
(1) into consideration, we have $s_{ij}=s_{0}$, $t_{ij}=t_{0}$, and
$z=\sum_{i=0}^{n-1}\sum_{j=0}^{m-1}u_{ij}\alpha_{ij}^{s_{0}}\beta_{ij}^{t_{0}}$
with
$u_{ij}=(-1)^{(i+j)(ms_{0}+nt_{0})}(\prod_{l=0}^{j-1}\prod_{p=0}^{n-1}q_{pl}^{t_{0}})
(\prod_{p=0}^{i-1}\prod_{l=0}^{m-1}q_{pl}^{-s_{0}})u_{00}$.
Moreover, Since $u_{00}\neq 0$, we have
$\eta^{t_0}=(-1)^{m(ms_0+nt_0)}$ and $
\eta^{s_0}=(-1)^{n(ms_0+nt_0)}.$

Finally, we consider the condition (3). If $|z|=0$, then
$u_{-1}c_{ij}=c_{ij}z=zc_{ij}=u_{00}c_{ij}$, which yields
$u_{00}=u_{-1}$, and thus $z=u_{00}\in k$. If  $|z|\ne 0$, then
$z=\sum_{i=0}^{n-1}\sum_{j=0}^{m-1}u_{ij}\alpha_{ij}^{s_{0}}\beta_{ij}^{t_{0}}$.
Thus $0=c_{ij}z=(-1)^{|z|}zc_{ij}$ forces $t_0\geq 1$ as desired
because  $\beta_{ij}c_{ij}$ lie in $\tilde{I}_q^{\perp}$ but
$\alpha_{ij}c_{ij}$ do not for $i\in \mathbb{Z}_n, j\in
\mathbb{Z}_m$. The proof of this lemma is finished. \hfill$\square$

\medskip

With a similar argument as in the proof of Proposition 2.2, if
$z\notin k$, then $\eta^{t_0}=(-1)^{m(ms_0+nt_0)}$ and
$\eta^{s_0}=(-1)^{n(ms_0+nt_0)}$,  which implies that $\eta$ is a
root of unity. Thus we immediately have

{\bf Proposition 4.2.}~If $\eta$ is not a root of unity, then
$Z_{gr}(E(\Gamma^{m,n}_q))=k$.

\medskip

{\bf Proposition 4.3.}~ Let
$\eta=\prod_{i=0}^{n-1}\prod_{j=0}^{m-1}q_{ij}$ be a primitive
$d$-th root of unity. If $\ch k=2$ or $m, n$ are even, then
$Z_{gr}(E(\Gamma^{m,n}_q))\cong k\oplus k[x, y]y.$

{\bf Proof.}~In the case that $\ch k=2$ or $m, n$ are even, we have
$\eta^{s_0}=\eta^{t_0}=1$, and thus $d|s_0, d|t_0$ since $\eta$ is a
primitive $d$-th root of unity. We assume $s_0=sd, t_0=td$ for some
integers $s\geq 0$ and $t\geq 1$ by $t_0=td\geq 1$.

Recall that for any homogeneous element $z\in
Z_{gr}(E(\Gamma^{m,n}_q))$, if $z\notin k$, then
$$
\begin{aligned}
z&=u_{00}\sum\limits_{i=0}^{n-1}\sum\limits_{j=0}^{m-1}
(-1)^{(i+j)(ms_{0}+nt_{0})}(\prod_{p=0}^{i-1}\prod_{l=0}^{m-1}q_{pl}^{-s_{0}})
(\prod_{l=0}^{j-1}\prod_{p=0}^{n-1}q_{pl}^{t_{0}})
\alpha_{ij}^{s_{0}}\beta_{ij}^{t_{0}}\\
&=u_{00}\sum\limits_{i=0}^{n-1}\sum\limits_{j=0}^{m-1}
\big((\prod_{p=0}^{i-1}\prod_{l=0}^{m-1}q_{pl})^{-1}\alpha_{ij}\big)^{s_{0}}
\big((\prod_{l=0}^{j-1}\prod_{p=0}^{n-1}q_{pl})\beta_{ij}\big)^{t_{0}}\\
&=u_{00}\big(\sum\limits_{i=0}^{n-1}\sum\limits_{j=0}^{m-1}
(\prod_{p=0}^{i-1}\prod_{l=0}^{m-1}q_{pl})^{-1}\alpha_{ij}\big)^{s_{0}}
\big(\sum\limits_{i=0}^{n-1}\sum\limits_{j=0}^{m-1}(\prod_{l=0}^{j-1}
\prod_{p=0}^{n-1}q_{pl})\beta_{ij}\big)^{t_{0}}
\end{aligned}
$$
$$
\begin{aligned}
 &=u_{00}\big(\sum\limits_{i=0}^{n-1}\sum\limits_{j=0}^{m-1}
(\prod_{p=0}^{i-1}\prod_{l=0}^{m-1}q_{pl})^{-1}\alpha_{ij}\big)^{sd}
\big(\sum\limits_{i=0}^{n-1}\sum\limits_{j=0}^{m-1}(\prod_{l=0}^{j-1}
\prod_{p=0}^{n-1}q_{pl})\beta_{ij}\big)^{td}.
\end{aligned} \eqno(4\mbox{-}1)
$$
Set $x=\big(\sum\limits_{i=0}^{n-1}\sum\limits_{j=0}^{m-1}
(\prod_{p=0}^{i-1}\prod_{l=0}^{m-1}q_{pl})^{-1}\alpha_{ij}\big)^{d}$
and
$y=\big(\sum\limits_{i=0}^{n-1}\sum\limits_{j=0}^{m-1}(\prod_{l=0}^{j-1}
\prod_{p=0}^{n-1}q_{pl})\beta_{ij}\big)^{d}$. Then we have
$yx=\eta^{d^2}xy=xy$ by $\alpha_{ij}\beta_{ij}=\eta
\beta_{ij}\alpha_{ij}$, and thus $z$ can be written as a scalar
multiple of $x^s y^t$ with $t\geq 1$. In addition, since
$\{x^{s}y^{l-s}\mid 0\leq s\leq l\}$ is a linearly independent set
for any fixed $l$, there is no additional homogeneous relation in
$Z_{gr}(E(\Gamma^{m,n}_q))$, and hence
$Z_{gr}(E(\Gamma^{m,n}_q))\cong k\oplus k[x, y]y.$
 \hfill$\square$

\medskip

{\bf Proposition 4.4.} Suppose that $\eta$ is a primitive $d$-th
root of unity, $\ch k\neq 2$ and $m, n$ have the different parity.
Then $Z_{gr}(E(\Gamma^{m,n}_q))\cong k\oplus k[x,y]y.$

{\bf Proof.} Without loss of generality, we assume that $n$ is even
and $m$ is odd. Then, by Lemma 4.1, the equalities
$\eta^{t_0}=(-1)^{m(ms_0+nt_0)}$ and $
\eta^{s_0}=(-1)^{n(ms_0+nt_0)}$ imply $\eta^{s_0}=1$ and
$\eta^{t_0}=(-1)^{ms_0}$. And thus we can write $s_0=sd$ for some
integer $s$.

(i) If $d$ is even, then $\eta^{t_0}=(-1)^{ms_0}=(-1)^{msd}=1$, thus
$d|t_0$ as well. With the same argument as that in the proof of
Proposition 4.3 we have $Z_{gr}(E(\Gamma^{m,n}_q))\cong k\oplus
k[x,y]y$ as desired.

(ii) If $d$ is odd, then $\eta^{2t_0}=1$, which implies that
$d|2t_0$ and thus $d|t_0$. We assume that $t_0=td$ with $t\geq 1$.
Since $1=\eta^{t_0}=(-1)^{ms_0}$ and $m$ is odd, we have $s_0$ is
even, and $s_0=sd$ implies that $s$ is even as well. As what we have
done in the proof of Proposition 4.3, for any homogeneous element
$z\in Z_{gr}(E(\Gamma^{m,n}_q))\setminus k$, we have the equality
$(4\mbox{-}1)$. Set
$x=\big(\sum\limits_{i=0}^{n-1}\sum\limits_{j=0}^{m-1}
(\prod_{p=0}^{i-1}\prod_{l=0}^{m-1}q_{pl})^{-1}\alpha_{ij}\big)^{2d}$
and $y=\big(\sum\limits_{i=0}^{n-1}\sum\limits_{j=0}^{m-1}
(\prod_{l=0}^{j-1}\prod_{p=0}^{n-1}q_{pl})\beta_{ij}\big)^{d}$. Then
we have $yx=\eta^{2d^2}xy=xy$ and $z=u_{00}x^{s/2} y^t$ with
$u_{00}\in k^*$,  $t\geq 1$ and $s/2=0,1,2,\cdots$. Again, there is
no additional homogeneous relation in $Z_{gr}(E(\Gamma^{m,n}_q))$.
Therefore, $Z_{gr}(E(\Gamma^{m,n}_q))\cong k\oplus k[x,y]y.$
\hfill$\square$

\medskip

{\bf Proposition 4.5.}~ Let $\eta$ be a primitive $d$-th root of
unity. If $\ch k\neq 2$ and both $m$ and $n$ are odd, then
$$Z_{gr}(E(\Gamma^{m,n}_q))\cong \left\{\begin{array}{ll}(k\oplus k[x,y]y)^{ev}, &\mbox{if}\; d \;\mbox{is odd};\\
k\oplus k[x,y]y, &\mbox{otherwise,}\end{array}\right.
 $$ where $(k\oplus k[x,y]y)^{ev}$ denotes the subalgebra of $k\oplus k[x,y]y$ spanned by all
even degree homogeneous elements as $k$-vector space.

{\bf Proof.}~If $\ch k\neq 2$ and $m, n$ are odd, then
$\eta^{s_0}=\eta^{t_0}=(-1)^{s_0+t_0}$, and thus $d|2s_0, d|2t_0$.

(i) In the case that $d$ is odd, we have $d|s_0$ and $d|t_0$. We
assume that $s_0=sd$ and $t_0=td$ for $s\geq 0$ and $t\geq 1$.
Moreover, $1=\eta^{s_0}=\eta^{t_0}=(-1)^{s_0+t_0}=(-1)^{s+t}$
implies that $s+t$ is even. In a similar way to the proof of
Proposition 4.3,  set
$x=\big(\sum\limits_{i=0}^{n-1}\sum\limits_{j=0}^{m-1}
(\prod_{p=0}^{i-1}\prod_{l=0}^{m-1}q_{pl})^{-1}\alpha_{ij}\big)^{d}$
and $y=\big(\sum\limits_{i=0}^{n-1}\sum\limits_{j=0}^{m-1}
(\prod_{l=0}^{j-1}\prod_{p=0}^{n-1}q_{pl})\beta_{ij}\big)^{d}$. Then
we have that $yx=\eta^{d^2}xy=xy$, and that $z\in
Z_{gr}(E(\Gamma^{m,n}_q))\setminus k$ can be written as a scalar
multiple of $x^s y^t$ with $t\geq 1$ and $s+t$ is even. Moreover,
$yx=xy$ is the sole relation in $Z_{gr}(E(\Gamma^{m,n}_q))$. So
$Z_{gr}(E(\Gamma^{m,n}_q))\cong (k\oplus k[x,y]y)^{ev}$.

(ii)~If $d$ is even, then $(d/2)|s_0$ and $(d/2)|t_0$. We write
$s_0=s(d/2)$ and $t_0=t(d/2)$ with $s\geq 0$ and $t\geq 1$. By
$\eta^{s_0}=\eta^{t_0}=(-1)^{s_0+t_0}$, we have $\eta^{s_0+t_0}=1$,
which implies that $d|{s_0+t_0}$. Since $s_0+t_0=d(s+t)/2$, we have
$s+t$ is even. Thus
$1=(-1)^{d(s+t)/2}=(-1)^{s_0+t_0}=\eta^{s_0}=\eta^{t_0}$, which
yields $d|s_0$ and $d|t_0$. Therefore, the rest of the proof in this
case is the same as the proof of the Proposition 4.3 and we omit it.
So $Z_{gr}(E(\Gamma^{m,n}_q))\cong k\oplus k[x,y]y$.
\hfill$\square$

\medskip

From the above four propositions, we have $\mathcal{N}_Z=0$, where
$\mathcal{N}_Z$ denotes the ideal of $Z_{gr}(E(\Gamma^{m,n}_q))$
generated by all nilpotent elements. By the isomorphism
 $\HH^{\ast}(\Gamma_q)/\mathcal{N}\cong
Z_{gr}(E(\Gamma_q))/\mathcal{N}_Z$ in \cite{SS, BGSS}, we have
$\HH^{\ast}(\Gamma_q)/\mathcal{N}\cong Z_{gr}(E(\Gamma_q))$.
Therefore, as  is shown in the following theorem,  $\Gamma_q^{m,n}$
provides more counterexamples to Snashall-Solberg's conjecture.

{\bf Theorem 4.6.}~ Let $\Gamma_q^{m,n}$ be the algebras defined in
the beginning of this section. Then

$$
\HH^{\ast}(\Gamma^{m,n}_q)/\mathcal{N}\cong\left\{\begin{array}{lll}
k,
&\mbox{if } \eta \mbox{ is not a root of unity;}\\[1ex]
(k\oplus k[x,y]y)^{ev}, &{\mbox{if}\; \ch k\neq 2, \;\eta \mbox{ is
a } d \mbox{-th primitive root of unity } \atop \mbox{and }d, m, n
\mbox{ is
odd};}\\[1ex]
k\oplus k[x,y]y, &\mbox{otherwise.}\end{array}\right.
 $$ As a consequence, if $\eta$ is a root of unity, then
$\HH^{\ast}(\Gamma_q^{m,n})/\mathcal{N}$ is not finitely generated
as algebra.

{\bf Proof.}~The first part of this theorem follows directly from
Propositions  4.2-4.5, and the proof of the second part is similar
to that of Theorem 2.4. \hfill$\square$

\medskip
{\bf Remark.} Our result is still true when $m=1$ or $n=1$.
Moreover, if $m=n=1$, the above result coincides with that of
\cite{Sna, XZ}.

\bigskip
\noindent{\bf\large Appendix. }
 \medskip

In this appendix we give a complete proof of Proposition 2.3, which
is a bit subtle modification of the proofs of the propositions 2.4
and 2.5 in \cite{ST}.

{\it Proof of Proposition 2.3.}~ We divide into two cases to finish
the proof.

 {\bf Case 1.} $m$ is even or $\ch
k=2$. In this case we have $\zeta^{t_0}=\zeta^{s_0}=1$. Since
$\zeta$ is a primitive $d$-th root of unity,  $d|s_0$ and $ d|t_0$.
We recall that $s_0\equiv t_0(\mod ~m)$, so $t_0=rm+s_0$, for some
integer $r$. Moreover, we have $u_1=(-1)^{s_0}(q_1\cdots
q_{t_0})^{-1}u_0=(-1)^{t_0}(q_1\cdots q_{s_0})^{-1}u_0$. If $m$ is
even or $\ch k=2$, then $(-1)^{s_0}=(-1)^{t_0}$, and thus
$q_1q_2\cdots q_{s_0}=q_1q_2\cdots q_{t_0}$. If $s_0\geq t_0$, then
$q_{t_0+1}q_{t_0+2}\cdots q_{s_0}=1$; on the other hand, if $t_0\geq
s_0$, then $q_{s_0+1}q_{s_0+2}\cdots q_{t_0}=1$. So in both cases,
we have $\zeta^r=1$, and thus $d|r$. So we can write $t_0=dhm+s_0$
for some integer $h$.

For any $z\in Z_{gr}(E(\Gamma_q))$, if $z$ is not in $k$,
$z=\sum_{i=0}^{m-1}u_i\gamma^{s_0}_i\delta^{t_0}_i$ with $t_0\geq
1$, and $u_i=(-1)^{is_0}\prod_{k=1}^{i}(q_k\cdots
q_{k+t_0-1})^{-1}u_0=(-1)^{it_0}\prod_{k=1}^{i}(q_k\cdots
q_{k+s_0-1})^{-1}u_0$ for $i=1,2,\dots, m-1$.

(i) ~We first consider the case $s_0=0$ and $t_0\geq 1$. Then
$z=\sum_{i=0}^{m-1}u_i\delta^{t_0}_i=\sum_{i=0}^{m-1}u_i\delta^{dhm}_i$
with $u_i=(-1)^{it_0}u_0=u_0$, and thus
$$z=\sum_{i=0}^{m-1}u_0\delta^{dhm}_i=u_0(\sum_{i=0}^{m-1}\delta^{dm}_i)^h.$$

(ii)~When $s_0, t_0\geq 1$, without loss of generality, we may
assume $s_0\leq t_0$. Then
\begin{eqnarray*}z&=& \sum_{i=0}^{m-1}(-1)^{it_0}\prod_{k=1}^{i}(q_k\cdots
q_{k+s_0-1})^{-1}u_0\gamma^{s_0}_i\delta^{t_0}_i\\
&=&\sum_{i=0}^{m-1}(-1)^{is_0}\prod_{k=1}^{i}(q_k\cdots
q_{k+s_0-1})^{-1}u_0\gamma^{s_0}_i\delta^{s_0}_i\delta^{dhm}_i\\
&=&u_0\Big(\sum_{i=0}^{m-1}(-1)^{is_0}\prod_{k=1}^{i}(q_k\cdots
q_{k+s_0-1})^{-1}\gamma^{s_0}_i\delta^{s_0}_i\Big)\Big(\sum_{i=0}^{m-1}\delta^{dm}_i\Big)^{h}.
\end{eqnarray*}
We assume $s_0=\alpha dm+s$, $0\leq s\leq dm-1$. Then
$(-1)^{s_0}=(-1)^s$, and $q_kq_{k+1}\cdots q_{k+s_0-1}=\zeta^{\alpha
d}q_kq_{k+1}\cdots q_{k+s-1}=q_kq_{k+1}\cdots q_{k+s-1}$. And the
above equality changes into
$$
\begin{array}{ccl}z
&=&u_0\Big(\sum_{i=0}^{m-1}(-1)^{is}\prod_{k=1}^{i}(q_k\cdots
q_{k+s-1})^{-1}\gamma^{\alpha dm+s}_i\delta^{\alpha dm+s}_i\Big)\Big(\sum_{i=0}^{m-1}\delta^{dm}_i\Big)^{h}\\
&=&u_0\Big(\sum_{i=0}^{m-1}(-1)^{is}\prod_{k=1}^{i}(q_k\cdots
q_{k+s-1})^{-1}\gamma^{s}_i
(\sum_{i=0}^{m-1}\gamma^{\alpha dm}_i\delta^{\alpha dm}_i)\delta^{s}_i\Big)\Big(\sum_{i=0}^{m-1}\delta^{dm}_i\Big)^{h}\\
&=&u_0\Big(\sum_{i=0}^{m-1}(-1)^{is}\prod_{k=1}^{i}(q_k\cdots
q_{k+s-1})^{-1}\gamma^{s}_i
\delta^{s}_i\Big)\Big(\sum_{i=0}^{m-1}\gamma^{dm}_i\Big)^{\alpha
}\Big(\sum_{i=0}^{m-1}\delta^{dm}_i\Big)^{\alpha+h}.
\end{array}\eqno{(A\mbox{-}1)}
$$
Since $d|s_0$ and $s_0=\alpha dm+s$, we have $d|s$ and $0\leq s\leq
dm-1$, and thus $s\in \{0, d, 2d, \cdots, d(m-1)\}$. We assume
$s=jd$,
 and define
$$z_j=\sum_{i=0}^{m-1}(-1)^{ijd}\prod_{k=1}^{i}(q_k\cdots
q_{k+jd-1})^{-1}\gamma^{jd}_i \delta^{jd}_i$$ for $0\leq j\leq m$.
In particular, $z_0=1$. Moreover, we have that
\begin{eqnarray*}z_jz_1&=&\Big(\sum_{i=0}^{m-1}(-1)^{ijd}\prod_{k=1}^{i}(q_k\cdots
q_{k+jd-1})^{-1}\gamma^{jd}_i
\delta^{jd}_i\Big)\Big(\sum_{i=0}^{m-1}(-1)^{id}\prod_{k=1}^{i}(q_k\cdots
q_{k+d-1})^{-1}\gamma^{d}_i \delta^{d}_i\Big)\\
&=&\sum_{i=0}^{m-1}(-1)^{i(j+1)d}\prod_{k=1}^{i}(q_k\cdots
q_{k+jd-1})^{-1}\prod_{k=1}^{i}(q_k\cdots
q_{k+d-1})^{-1}\gamma^{jd}_i \delta^{jd}_i\gamma^{d}_i \delta^{d}_i\\
 &=&(-1)^{jd}\prod_{k=1}^{jd}(q_k\cdots
q_{k+d-1})^{-1}\Big(\sum_{i=0}^{m-1}(-1)^{i(j+1)d}\prod_{k=1}^{i}(q_k\cdots
q_{k+(j+1)d-1})^{-1}\gamma^{(j+1)d}_i
\delta^{(j+1)d}_i\Big)\\
&=&(-1)^{jd}\prod_{k=1}^{jd}(q_k\cdots q_{k+d-1})^{-1}z_{j+1}.
\end{eqnarray*} Thus we have $z_1^j=(-1)^{\sum_{i=1}^{j-1}id}(\prod_{l=1}^{j-1}\prod_{k=1}^{ld}(q_k\cdots
q_{k+d-1})^{-1})z_j$, for $j=1,2,\cdots,m$. In particular,
\begin{eqnarray*}z_1^m&=&(-1)^{\sum_{i=1}^{m-1}id}\Big(\prod_{l=1}^{m-1}\prod_{k=1}^{ld}(q_k\cdots
q_{k+d-1})^{-1}\Big)z_m\\
&=&(-1)^{md/2}\Big(\prod_{l=1}^{m-1}\prod_{k=1}^{ld}(q_k\cdots
q_{k+d-1})^{-1}\Big)\Big(\sum_{i=0}^{m-1}(-1)^{imd}\prod_{k=1}^{i}(q_k\cdots
q_{k+md-1})^{-1}\gamma^{md}_i \delta^{md}_i\Big)\\
&=&(-1)^{md/2}\Big(\prod_{l=1}^{m-1}\prod_{k=1}^{ld}(q_k\cdots
q_{k+d-1})^{-1}\Big)\Big(\sum_{i=0}^{m-1}\gamma^{md}_i\Big)\Big(\sum_{i=0}^{m-1}\delta^{md}_i\Big).
\end{eqnarray*}
Set $x=\sum_{i=0}^{m-1}\gamma^{md}_i$,
$y=\sum_{i=0}^{m-1}\delta^{md}_i$,
$w=z_1=\sum_{i=0}^{m-1}(-1)^{id}\prod_{k=1}^{i}(q_k\cdots
q_{k+d-1})^{-1}\gamma^{d}_i \delta^{d}_i$ and
$\epsilon_d=(-1)^{md/2}\prod_{l=1}^{m-1}\prod_{k=1}^{ld}(q_k\cdots
q_{k+d-1})^{-1}$. Then $w^m=\epsilon_d xy$. Moreover, by the formula
(A-1), we have  $z\in k$ or $z$ has the form
$z=u'_0w^jx^{\alpha}y^{\alpha+h}$for any homogeneous element $z\in
Z_{gr}(E(\Gamma_q))$, where $u'_0\in k$, $s_0=s+\alpha dm=(j+\alpha
m)d>0$, that is, $j+\alpha m>0$, and thus $j+\alpha>0$. Similarly,
if $s_0\geq t_0$, then $z=u''_0w^jx^{\alpha+h}y^{\alpha}$ wit
h$u''_0\in k$ and $j+\alpha>0$. Therefore, in both cases, any
homogeneous element $z\in Z_{gr}(E(\Gamma_q))\setminus k$ can be
written as a scalar multiple of $x^iy^jw^l$ with $j+l>0$ and
$w^m=\epsilon_d xy$. In particular, any scalar multiple of $x^i$
does not lie in $Z_{gr}(E(\Gamma_q))$, for $i=1, 2, \cdots$.

As what Snashall and Taillefer have done in \cite[Lemma 2.3]{ST}, we
 claim that the elements $x, y, w$ don't have additional
relation except $w^m=\epsilon_d xy$ in $Z_{gr}(E(\Gamma_q))$.

Indeed, since the elements $x^iy^{n-i}$ have different degree, for
$i=0, 1, \cdots, n $, thus they are linearly independent in
$Z_{gr}(E(\Gamma_q))$. So any additional relation in
$Z_{gr}(E(\Gamma_q))$ is length homogeneous of the form
$$f_0(x,y)+f_1(x,y)w+\cdots +f_{m-1}w^{m-1}=0, \eqno{(A\mbox{-}2)}$$
where $f_i(x,y)=\sum_{j=0}^{n_i}k_{ij}x^jy^{n_i-j}\in k[x, y]$, and
$|f_0(x,y)|=|f_1(x,y)|+|w|$, which implies $n_0|y|=n_1|y|+|w|$, and
thus $n_0md=n_1md+2d$, that is, $n_0m=n_1m+2$.

If $m=1$, then $w=\epsilon_d xy$, and thus any element $z\in
Z_{gr}(E(\Gamma_q))$ can be generated by $x, y$. So there is no
additional relation in $Z_{gr}(E(\Gamma_q))$.

Now we consider the case $m\geq 2$. $n_0m=n_1m+2$ implies $m=2$ and
$n_0=n_1+1$. Then $|x|=|y|=|w|=2d$, and we may choose the minimal
$n_0$ such that $f_0(x, y)+f_1(x, y)w=0$ with $|f_0(x, y)|=2n_0d$
and $|f_1(x,y)|=2(n_0-1)d$. Since $x^{n_0} \notin
Z_{gr}(E(\Gamma_q))$, $f_0(x,
y)=\sum_{j=0}^{n_0-1}k_{0j}x^jy^{n_0-j}$ and $f_1(x,
y)=\sum_{i=0}^{n_0-1}k_{1j}x^jy^{n_0-j-1}$. Then
$f^2_0(x,y)=f^2_1(x,y)w^2=\epsilon_d f^2_1(x, y)xy$. Comparing the
coefficients of  $y^{2n_0}$ and $x^{2n_0-1}y$, we have
$k_{00}=k_{1,n_0-1}=0$, and then $f_1(x,
y)=\sum_{j=0}^{n_0-2}k_{1j}x^jy^{n_0-j-1}$ and $f_0(x,
y)=\epsilon_d^{-1}f_0^{'}(x,y)w^2$ with
$f_0^{'}(x,y)=\sum_{j=0}^{n_0-1}k_{0,j+1}x^jy^{n_0-j-1}$, thus
$\epsilon_d^{-1}f_0^{'}(x,y)w+f_1(x, y)=0$, which contradicts to the
minimality of $n$.

{\bf Case 2.} $m$ is odd and char$k\ne 2$.  By the conditions
$\zeta^{s_0}=(-1)^{mt_0}$ and $\zeta^{t_0}=(-1)^{ms_0}$, we know
that $\zeta^{2s_0}=\zeta^{2t_0}=1$. Since $\zeta$ is a primitive
$d$-th root of unity, $d|2s_0$, and $d|2t_0$. Recall that $s_0\equiv
t_0(\mod ~ m)$, that is,  $s_0=t_0+rm$ for some integer $r$, and
$u_1=(-1)^{s_0}(q_1\cdots q_{t_0})^{-1}u_0=(-1)^{t_0}(q_1\cdots
q_{s_0})^{-1}u_0$. If $s_0\geq t_0$, then $q_{t_0+1}q_{t_0+2}\cdots
q_{s_0}=(-1)^{s_0-t_0}$; on the other hand, if $t_0\geq s_0$, then
$q_{s_0+1}q_{s_0+2}\cdots q_{t_0}=(-1)^{t_0-s_0}$. So, in both
cases, we have $\zeta^{2r}=1$, and thus $d|2r$. Then
$dm|2(t_0-s_0)$. We assume that $s_0=\alpha dm+s$ and $t_0=\beta
dm+t$, where $0\leq s, t \leq dm-1$, then $dm|2(s-t)$, without loss
of generality, we assume $s\geq t$, then $2(s-t)=0$ or $2(s-t)=dm$.

Now we assert that $2(s-t)=0$ and thus $t=s$. Otherwise, we will
have $2(s-t)=dm$. Since $m$ is odd and $d$ is even,   $s-t$ and
$s_0-t_0$ have the same parity. Moreover,
$(-1)^{s_0-t_0}=\zeta^r=\zeta^{(s_0-t_0)/m}=\zeta^{(\alpha-\beta)d+(s-t)/m}
=\zeta^{(s-t)/m}=\zeta^{d/2}=-1$. Therefore, $s-t$ is odd and $d/2$
is odd. We can also get the equality
$(-1)^{s_0+t_0}=(-1)^{m(s_0+t_0)}=(-1)^{ms_0}(-1)^{mt_0}=\zeta^{s_0+t_0}=\zeta^{s+t}=
\zeta^{2t+s-t}=\zeta^{2t+(dm)/2}=\zeta^{2t}(-1)^m=-\zeta^{2t}$. So
$\zeta^{4t}=1$, and thus $d|4t$ and $(d/2)|2t$. Moreover, since
$d/2$ is odd, $(d/2)|t$. We assume that $t=ld/2$ for some integer
$l$. If $t$ is even, then $l$ is even, and we have
$1=(-1)^t=(-1)^{t_0}=\zeta^{t_0}=\zeta^{s}=\zeta^{t+(s-t)}=\zeta^{(l+m)d/2}=\zeta^{l+m}=-1$,
this yields a contradiction. Therefore, $t$ is odd, then $l$ is odd,
$s=t+(s-t)=(l+m)d/2$ is even, and we have
$1=(-1)^s=(-1)^{s_0}=\zeta^{t_0}=\zeta^t=\zeta^{ld/2}=(-1)^l=-1$, a
contradiction again. So $2(s-t)=0$ and thus $t=s$ as desired.

Since $t_0=\alpha dm+t, s_0=\beta dm+s$ and $t=s$, we have
$s_0-t_0=(\alpha-\beta)dm$ and
$1=\zeta^{(\alpha-\beta)dm}=\zeta^{s_0-t_0}=(-1)^{m(t_0-s_0)}
=(-1)^{t_0-s_0}=(-1)^{(\beta-\alpha)dm}$. So $\alpha dm$ and $\beta
dm$ have the same parity, and thus $\alpha d$ and $\beta d$ have the
same parity. By squaring the equality
$\zeta^t=\zeta^{t_0}=(-1)^{ms_0}$, we know $\zeta^{2t}=1$, and thus
$d|2t$ with $0\leq 2t<2dm$. We assume $2t=dl$ for some integer
$0\leq l<2m$.

Now, we will describe any homogeneous element in
$Z_{gr}(E(\Gamma_q))$. We recall that if $z$ is not in $k$, $z=
\sum_{i=0}^{m-1}(-1)^{it_0}\prod_{k=1}^{i}(q_k\cdots
q_{k+s_0-1})^{-1}u_0\gamma^{s_0}_i\delta^{t_0}_i $ with $t_0\geq 1$.

(i)~If $d$ is odd, then by $2t=dl$, we have $l$ is even and since
$\alpha d, \beta d$ have the same parity,
$1=\zeta^{dl/2}=\zeta^{t}=\zeta^{t_0}=(-1)^{ms_0}=(-1)^{s_0}=(-1)^{\alpha
dm+s}=(-1)^{\alpha d+t}=(-1)^{(\alpha+l/2)d}=(-1)^{(\beta+l/2)d}$.
So $(\alpha+l/2)d$ and  $(\beta+l/2)d$ are even with $0\leq l<m$. If
$\alpha $ is even, then $l/2$ and thus $t_0=\beta dm+dl/2$ is even.
So we have
\begin{eqnarray*}z&=&
u_0\sum_{i=0}^{m-1}\prod_{k=1}^{i}(q_k\cdots
q_{k+s_0-1})^{-1}\gamma^{s_0}_i\delta^{t_0}_i\\
&=&u_0\sum_{i=0}^{m-1}\prod_{k=1}^{i}(q_k\cdots
q_{k+dl/2-1})^{-1}\gamma^{\alpha dm+dl/2}_i\delta^{\beta dm+dl/2}_i\\
&=&u_0\sum_{i=0}^{m-1}\prod_{k=1}^{i}(q_k\cdots
q_{k+dl/2-1})^{-1}\gamma^{dl/2}_i\delta^{dl/2}_i(\sum_{i=0}^{m-1}\gamma^{2dm}_i)^{\alpha/2}
(\sum_{i=0}^{m-1}\delta^{2dm}_i)^{\beta/2}.
\end{eqnarray*}
Similarly, if $\alpha$ is odd, then $l/2$ is odd, and $t_0$ is even,
we have
\begin{eqnarray*}z&=&
u_0\sum_{i=0}^{m-1}\prod_{k=1}^{i}(q_k\cdots
q_{k+s_0-1})^{-1}\gamma^{s_0}_i\delta^{t_0}_i\\
&=&u_0\sum_{i=0}^{m-1}\prod_{k=1}^{i}(q_k\cdots
q_{k+dl/2-1})^{-1}\gamma^{\alpha dm+dl/2}_i\delta^{\beta dm+dl/2}_i
\end{eqnarray*}
\begin{eqnarray*}
&=&\sum_{i=0}^{m-1}\prod_{k=1}^{i}(q_k\cdots
q_{k+dl/2-1})^{-1}u_0\gamma^{d(l/2+m)}_i\delta^{d(l/2+m)}_i(\sum_{i=0}^{m-1}\gamma^{2dm}_i)^{(\alpha-1)/2}
(\sum_{i=0}^{m-1}\delta^{2dm}_i)^{(\beta-1)/2}.
\end{eqnarray*}
As what we have done in the case 1, we define
$$z_j=\sum_{i=0}^{m-1}\prod_{k=1}^{i}(q_k\cdots
q_{k+2dj-1})^{-1}\gamma^{2dj}_i\delta^{2dj}_i$$ for $j=1, 2, \cdots,
m$, then $z_1=\sum_{i=0}^{m-1}\prod_{k=1}^{i}(q_k\cdots
q_{k+2d-1})^{-1}\gamma^{2d}_i\delta^{2d}_i$. Moreover, by a
straightforward verification, $$z_1z_j=\prod_{k=1}^{2dj}(q_k\cdots
q_{k+2d-1})^{-1}z_{j+1},$$ for $j=1, 2, \cdots, m$. Thus,
$z_{1}^j=\prod_{l=1}^{j-1}\prod_{k=1}^{2dl}(q_k\cdots
q_{k+2d-1})^{-1}z_j$, for $j=1, 2, \cdots, m$. In particular,
\begin{eqnarray*}z_1^m&=&\prod_{l=1}^{m-1}\prod_{k=1}^{2dl}(q_k\cdots
q_{k+2d-1})^{-1}z_m\\
&=&\prod_{l=1}^{m-1}\prod_{k=1}^{2dl}(q_k\cdots
q_{k+2d-1})^{-1}(\sum_{i=0}^{m-1}\gamma^{2dm}_i)
(\sum_{i=0}^{m-1}\delta^{2dm}_i). \end{eqnarray*}

 Set
$x=\sum_{i=0}^{m-1}\gamma^{2md}_i$,
$y=\sum_{i=0}^{m-1}\delta^{2md}_i$, $w=z_1$ and
$\epsilon_d=\prod_{l=1}^{m-1}\prod_{k=1}^{2dl}(q_k\cdots
q_{k+2d-1})^{-1}$. Then $w^m=\epsilon_d xy$. Moreover, if $\alpha $
is even, then any $z\in Z_{gr}(E(\Gamma_q))\setminus k$  is a scalar
multiple of $x^{\alpha/2}y^{\beta/2}w^{l/4}$ with $\beta/2+l/4>0$
(because $t_0=\beta dm+dl/2>0$). Similarly, if $\alpha$ is odd and
$z\in Z_{gr}(E(\Gamma_q))\setminus k$, then $z$ is a scalar multiple
of $x^{(\alpha-1)/2}y^{(\beta-1)/2}w^{(l/2+m)/2}$ with
$(\beta-1)/2+(l/2+m)/2>0$. In both cases, $z\in
Z_{gr}(E(\Gamma_q))\setminus k$  can be written as a scalar multiple
of $x^{i}y^{j}w^l$ with $j+l>0$ and $w^m=\epsilon_d xy$. Note that
any scalar multiple of $x^i$ does not belong to
$Z_{gr}(E(\Gamma_q))$, for $i=1, 2, \cdots$.

 With a similar argument as in the case 1, we can
assert that $x, y, w$ have no additional homogeneous relation except
$w^m=\epsilon_d xy$. Indeed, it suffices to note that
$n_0|x|=2n_0dm=|f_0(x,y)|=|f_1(x,y)|+|w|=n_1|y|+|w|=2n_1dm+4d$, and
thus $(n_0-n_1)m=2$ has no
 solution in $\mathbb{Z}$. So there is no additional homogeneous
 relation of the form (A-2) as required.

 (ii)~Now we consider the case  $d\equiv 0(\mod ~4)$. We assert that $l$ is even with $0\leq l/2<m$, and thus $t_0$ is
even. Otherwise, if $l$ is odd, then, by
$-1=(-1)^l=(\zeta^{(d/2)})^l=\zeta^{t}=\zeta^{t_0}=(-1)^{ms_0}=(-1)^{s_0}=(-1)^{\alpha
dm+ s}=(-1)^s=(-1)^t=(-1)^{dl/2}=(-1)^{d/2}$, we have that $d/2$ is
odd, which contradicts to $d\equiv 0(\mod ~4)$. Therefore, for any
given homogeneous element $z\in Z_{gr}(E(\Gamma_q))\setminus k$, we
have
\begin{eqnarray*}z
&=&u_0\sum_{i=0}^{m-1}\prod_{k=1}^{i}(q_k\cdots
q_{k+dl/2-1})^{-1}\gamma^{\alpha dm+dl/2}_i\delta^{\beta dm+dl/2}_i\\
&=&u_0\sum_{i=0}^{m-1}\prod_{k=1}^{i}(q_k\cdots
q_{k+dl/2-1})^{-1}\gamma^{dl/2}_i\delta^{dl/2}_i(\sum_{i=0}^{m-1}\gamma^{dm}_i)^{\alpha}
(\sum_{i=0}^{m-1}\delta^{dm}_i)^{\beta}.
\end{eqnarray*}
We define $z_j=\sum_{i=0}^{m-1}\prod_{k=1}^{i}(q_k\cdots
q_{k+dj-1})^{-1}\gamma^{dj}_i\delta^{dj}_i$, for $j=1, 2, \cdots,
m$. Then it is clear that $z_jz_1=\prod_{k=1}^{dj}(q_k\cdots
q_{k+d-1})^{-1}z_{j+1},$ for $j=1, 2, \cdots, m$. Thus,
$z_{1}^j=\prod_{l=1}^{j-1}\prod_{k=1}^{dl}(q_k\cdots
q_{k+d-1})^{-1}z_j$, for $j=1, 2, \cdots, m$.

Set $x=\sum_{i=0}^{m-1}\gamma^{md}_i$,
$y=\sum_{i=0}^{m-1}\delta^{md}_i$, $w=z_1$ and
$\epsilon_d=\prod_{l=1}^{m-1}\prod_{k=1}^{dl}(q_k\cdots
q_{k+d-1})^{-1}$. Then $w^m=\epsilon_d xy$. And we can write any
homogeneous element $z\in Z_{gr}(E(\Gamma_q))\setminus k$ as a
scalar multiple of $x^{\alpha}y^{\beta}w^{l/2}$ with $\beta +l/2>0$.
In particular, any scalar multiple of $x^i$ is not in
$Z_{gr}(E(\Gamma_q))$, for $i=1, 2, \cdots$.

Similarly, we can also prove that $x, y, w$ have no additional
relation except $w^m=\epsilon_d xy$. Thus we have
$Z_{gr}(E(\Gamma_q))\cong (k[x, y, w]/\langle w^m-\epsilon_d
xy\rangle)_{x^*},$ where $\epsilon_d=
\prod_{l=1}^{m-1}\prod_{k=1}^{ld}(q_k\cdots q_{k+d-1})^{-1}$.

(iii)~ If $d$ is even with $d\equiv 2(\mod ~ 4)$, then $d/2$ is odd,
and $t_0$ and $l$ have the same parity by $t_0=\beta dm + ld/2$,
where $0\leq l<2m$. So we can write any homogeneous element $z\in
Z_{gr}(E(\Gamma_q))$ that is not in $k$ as
\begin{eqnarray*}z
&=&u_0\sum_{i=0}^{m-1}\prod_{k=1}^{i}(-1)^{il}(q_k\cdots
q_{k+dl/2-1})^{-1}\gamma^{\alpha dm+dl/2}_i\delta^{\beta dm+dl/2}_i\\
&=&u_0\sum_{i=0}^{m-1}\prod_{k=1}^{i}(-1)^{il}(q_k\cdots
q_{k+dl/2-1})^{-1}\gamma^{dl/2}_i\delta^{dl/2}_i(\sum_{i=0}^{m-1}\gamma^{dm}_i)^{\alpha}
(\sum_{i=0}^{m-1}\delta^{dm}_i)^{\beta}.
\end{eqnarray*}
Similarly, define
$$z_j=\sum_{i=0}^{m-1}\prod_{k=1}^{i}(-1)^{ij}(q_k\cdots
q_{k+dj/2-1})^{-1}\gamma^{dj/2}_i\delta^{dj/2}_i,$$ for $j=1, 2,
\cdots, 2m$. Then we can verify that
$z_jz_1=\prod_{k=1}^{dj/2}(q_k\cdots q_{k+d/2-1})^{-1}z_{j+1},$ and
thus $z_{1}^j=\prod_{l=1}^{j-1}\prod_{k=1}^{dl/2}(q_k\cdots
q_{k+d/2-1})^{-1}z_j$, for $j=1, 2, \cdots, 2m$.

Set $x=\sum_{i=0}^{m-1}\gamma^{md}_i$,
$y=\sum_{i=0}^{m-1}\delta^{md}_i$, $w=z_1$ and
$\epsilon_d=\prod_{l=1}^{2m-1}\prod_{k=1}^{dl/2}(q_k\cdots
q_{k+d/2-1})^{-1}$. Then $w^{2m}=\epsilon_d xy$. And we can write
any  homogeneous element $z\in Z_{gr}(E(\Gamma_q))\setminus k$ as a
scalar multiple of $x^{\alpha}y^{\beta}w^{l}$ with $\beta+l>0$ since
$t_0=\beta dm +ld/2>0$.

Again,  there is no additional relation in $Z_{gr}(E(\Gamma_q))$
except $w^{2m}=\epsilon_d xy$, and  any scalar multiple of $x^i$ is
not in $Z_{gr}(E(\Gamma_q))$, for $i=1, 2, \cdots$.
 So we have
$Z_{gr}(E(\Gamma_q))\cong (k[x, y, w]/\langle w^{2m}-\epsilon_d
xy\rangle)_{x^*},$ where $\epsilon_d=
\prod_{l=1}^{2m-1}\prod_{k=1}^{ld/2}(q_k\cdots q_{k+d/2-1})^{-1}$ in
this case.\hfill$\square$

\bigskip
{\bf Acknowledgements}
\medskip

This research work was supported by the Natural Science Foundation
of China (Grant No. 10971206 and 11171325 ).


\medskip
\footnotesize





\begin{thebibliography}{99}

\bibitem{CE} H. Cartan,  S. Eilenberg, Homological Algebra,
 Princeton University Press, Princeton(1956)


\bibitem{Ger} M. Gerstenhaber, The cohomology struncture of an
associative ring, Ann. Math.  78(2) (1963) 267-288


\bibitem{Hap}D. Happel, Hochschild cohomology of finite-dimensional
algebras, Lecture Notes in Math. Springer, Berlin, 1404{1989}
108-126

\bibitem{BGMS} R.O. Buchweitz , E. L. Green, D. Madsen, {\O}. Solberg,
 Finite Hochschild cohomology without finite global dimension,
Math. Res. Lett. 12(2005) 805-816

\bibitem{BE}  P. A. Bergh, K. Erdmann, Homology and cohomology of quantum complete
intersections, Algebra and Number theory 2(5) (2008), 501-522

\bibitem{Opp} S. Oppermann,
Hochschild cohomology and homology of quantum complete
intersections, Algebra Number Theory 4(7)(2010)821-838

\bibitem{Han}Y. Han, Hochschild (co)homology dimension, J. London Math.
Soc.  73(2006)657-668

\bibitem{Car} J. F. Carlson,  Varieties and the cohomology ring of a
module, J. Algebra 5(1983)441-454

\bibitem{SS}  N. Snashall, {\O}. Solberg,  Support varieties and
Hochschild cohomology rings, Proc. London Math. Soc.  88(2004)
705-732


\bibitem{GS} E. L. Green, N. Snashall,  The Hochschild
cohomology ring modulo nipotence of a stacked monomial algebra,
Colloq. Math. 105(2006)233-258

\bibitem{GSS1}E. L. Green,  N. Snashall, {\O}. Solberg, The Hochschild
cohomology ring modulo nipotence of a monomial algebra,  J. Algebra
Appl.  5(2006)153-192

\bibitem{GSS2}E. L. Green, N. Snashall,  {\O}. Solberg, The Hochschild
cohomology ring of a selfinjective algebra of a finite
representation type, Proc. Amer. Math. Soc. 131(2003)3387-3393

\bibitem{Eve}L. Evens, The cohomology ring of a finite group,
Trans. Amer. Math. Soc. 101(1961)224-239

\bibitem{Ven} B. B. Venkov, Cohomology algebras for some classifying
space, Dokl. Akad. Nauk SSSR 127(1959) 943-944

\bibitem{XuF}F. Xu, Hochschild and ordinary cohomology rings of small
categories, Adv. Math.  219(2008)1872-1893

\bibitem{Sna}N. Snashall, Support varieties and the Hochschild
cohomology ring modulo nilpotence, Proceedings of the 41st Symposium
on Ring Theory and Representation Theory (2009)68-82

\bibitem{XZ} Y. Xu,  C. Zhang, Hochschild cohomology
 of a class of quantized Koszul algebras,
Preprint

\bibitem{ST0} N. Snashall, R. Taillefer, The Hochschild cohomology ring of 
 a class of special biserial algebras, J. Algebra Appl. 9(2010)73-122

\bibitem{ST} N. Snashall, R. Taillefer, Hochschild cohomology of socle
 deformations of a class of Koszul self-injective algebras,
 Colloquium Mathematicum, 119(2010)79-93

\bibitem{SP} A. Parker,  N. Snashall,
A family of Koszul self-injective algebras with finite Hochschild
cohomology. arXiv:1105.2215 [math.RA]

\bibitem{HZ}  D. Zhao, Y. Han, Koszul algebras and finite Galois
coverings, Sci. in China (Series A) 52(10)(2009)2145-2153

\bibitem{MMM1} E. N. Marcos, R. Mart\'inez-Villa, Ma. I. R. Martins, Hochschild
cohomology of skew group rings and invariants, Cent. Eur. J. Math.
2(2)(2004)177-190(electronic).

\bibitem{MMM2} E.N. Marcos, R. Mart¨ªnez-Villa, Ma.I.R. Martins, Addendum to:
Hochschild cohomology of skew group rings and invariants Cent. Eur.
J. Math. 2 (2) (2004) 177¨C190 (electronic), Cent. Eur. J. Math. 2
(4) (2004) 614 (electronic).

\bibitem{XT} Y. Xu  and X. Tang, Hochschild (co)homology of Galois coverings of
Grassmann algebras. Acta Math. Sinica, 25(10) (2009) 1693-1702

\bibitem{GHS} E. L. Green, J. R. Hunton, N. Snashall,
Coverings, the graded center and Hochschild cohomology, J.  Pure
Appl. Algebra 212(12)(2008)2691-2706

\bibitem{DMR}P. Dr\"axler, G. O. Michler, C. M. Ringel,
Computational methods for representations of groups and algebras,
Birkh\"auser Verlag, Berlin(1999)

\bibitem{GH}E. Green,  R. Q. Huang, Projective resolution of
straightening closed algebras generated by minors, Adv. Math.
110(1995) 314-333

\bibitem{BGSS}  R. Buchweitz, E. L. Green,  N. Snashall, {\O}. Solberg,
 Multiplicative structures for Koszul algebras.  Quart. J. Math.
59 (2008)441-454

\bibitem{Cib90} C. Cibils,  Rigidity of truncated quiver algebras,
Adv. Math.  79 (1990)18-42

\bibitem{HX}B. Hou, Y. Xu, Hochschild (co)homology of  $\mathbb{Z}_2\times \mathbb{Z}_2$-Galois
coverings of quantum exterior algebras, Bulletin of the Australian
Mathematican Society  78(1)(2008) 35-54

\bibitem{BK} M. C. R. Bulter, A. D. King, Minimal resolution of
algebras, J. algebra 212(1999)323-362

\bibitem{GHMS} E. L. Green, G. Hartman, E. N. Marcos,
{\O}. Solberg, Resolutions over Koszul algebras, Arch. Math.
85(2005)118-127





\end{thebibliography}
\end{document}